\documentclass[12pt]{article}

\usepackage{amsfonts,amsmath,graphicx}
\usepackage[final]{epsfig}
\usepackage{graphics}
\usepackage{amsmath}
\usepackage{amsfonts}
\usepackage{latexsym}
\usepackage{amssymb}
\usepackage{graphicx}
\usepackage{epstopdf}
\def\bb{\mathbb}
\def\frak{\mathfrak}
\def\scr{\mathcal}

\def\Ker{\mathop{\rm Ker}}
\def\Im{\mathop{\rm Im}}
\def\lk{\mathop{\rm lk}}
\def\ind{\mathop{\rm ind}\nolimits}

\def\span{\mathop{\rm span}}
\def\crit{\mathop{\rm crit}}
\def\R{\mathbb{R}}
\def\d{\partial}

\newtheorem{lemma}{Lemma}[section]

\newtheorem{theorem}{Theorem}[section]
\newtheorem{definition}{Definition}[section]
\newtheorem{corollary}{Corollary}[section]
\newtheorem{proposition}[theorem]{Proposition}

\title{Generating Families \\ and Legendrian Contact Homology in \\ the Standard Contact Space}

\author{Dmitry FUCHS\\University of California\\Davis
\and Dan RUTHERFORD \\Duke University\\}

\begin{document}

\date{}

\maketitle

\begin{abstract}  We show that if a Legendrian knot in standard contact ${\bb R}^3$ possesses a generating family then there exists an augmentation of the Chekanov-Eliashberg DGA so that the associated linearized contact homology (LCH) is isomorphic to singular homology groups arising from the generating family.  In this setting we show Sabloff's duality result for LCH may be viewed as Alexander duality. In addition, we provide an explicit construction of a generating family for a front diagram with graded normal ruling and give a new approach to augmentation $\Rightarrow$ normal ruling.

\end{abstract}

\tableofcontents

\section{Introduction}

Our main objects of study are Legendrian knots in the `standard contact space', $({\bb R}^3, \ker(dz- y dx ))$, and their invariants of Legendrian isotopy.  Aside from the underlying topological knot type there are two ``classical'' integer valued invariants known as the Thurston-Bennequin number and the rotation (or Maslov) number.  For some time it was not known whether there could exist distinct Legendrian knot types in the standard contact space with the same classical invariants (including topological knot type).  This question was answered in the affirmative independently by Chekanov and Eliashberg.  Motivated by Floer theory, they developed a rich algebraic invariant which takes the form of a DGA (Differential Graded Algebra) generated by self intersections of the knot's Lagrangian projection.   The differential is defined by counting immersed (or holomorphic) disks.  Generating families have provided a second source of non-classical invariants in the work of Traynor and her collaborators \cite{Tr} \cite{NgTr} \cite{JTr} and Chekanov-Pushkar \cite{ChP}.  The main results of this paper give an interpretation of a linearized version of the Chekanov and Eliashberg invariant in terms of generating families.  This strengthens previously known links between the two classes of invariants \cite{F} \cite{FI} \cite{K} \cite{NgTr} \cite{S} \cite{JTr}.

The Chekanov-Eliashberg DGA of a Legendrian knot $\ell$ and even its homology may be infinite dimensional.  An {\it augmentation} $\varepsilon$ is an algebra homomorphism of the Chekanov-Eliashberg DGA into ${\bb Z}_2$ which allows us to linearize the DGA to a finite dimensional complex with homology groups $H^\varepsilon_{*}(\ell)$.  These linearized homology groups may depend on the choice of $\epsilon$.  However, the set of graded groups $\{ H^\varepsilon_*(\ell) \}$ where $\varepsilon$ is any augmentation of the DGA provides a Legendrian isotopy invariant \cite{Ch}.  

Given a $1$-parameter family of functions $F: {\bb R}^n\times {\bb R} \rightarrow {\bb R}, f_t = F(\cdot, t)$ the fiber-wise critical set is immersed (under a transversality assumption) into the standard contact space as a Legendrian submanifold, $\ell$, and $F$ is called a {\it generating family} for $\ell$.  Traynor \cite{Tr} introduced a homological invariant obtained from generating families for a class of 2-component links in the solid torus and (with Jordan \cite{JTr}) in ${\bb R}^3$.  A variation\footnote{The authors became aware of this version of the generating family invariant through a letter from Peter Pushkar to the first author.   \cite{P}.  He considers the homology groups as a possibility for defining a Legendrian homology invariant, but as far as we know has never published anything in this regard.  The letter also suggests as a method for computing these groups a complex along the same lines as the one described in Section 5.} of the Traynor-Jordan invariant which applies without additional assumptions on the knot type is the following:

Given a generating family for $\ell$ consider the {\it difference function} \[\begin{array}{lr} w: {\bb R}^n\times{\bb R}^n\times{\bb R} \rightarrow {\bb R}, & w(x,y,t) = F(x,t)-F(y,t). \end{array}\]  Take $\delta$ close enough to $0$ so that the interval $(0, \delta)$ consists entirely of regular values of $w$ (0 itself is always a critical value of $w$) and consider the grading shifted homology groups ${\scr G}H_*(F) =H_{*+n+1}(w\geq \delta, w=\delta; {\bb Z}_2)$.  While the homology may depend on the choice of generating family, the set of graded groups $\{ {\scr G}H_*(F) \}$ where $F$ is any generating family (with restrictions on its behavior outside of  a compact set) for $\ell$ forms a Legendrian isotopy invariant.    

A central result of this paper is

\bigskip

\noindent {\bf Theorem 5.3}  {\it If $F$ is a linear at infinity generating family for $\ell$ then there exists a graded augmentation $\varepsilon$ for the Chekanov-Eliashberg DGA of $\ell$ such that ${\scr G}H_*(F) \cong H^\varepsilon_{*}(\ell)$.  That is,
$\{ {\scr G}H_*(F) \}\subset\{ H^\varepsilon_*(\ell) \} $.} 

\bigskip

It should be noted that a close variant of Theorem 5.3 is proved  in \cite{NgTr} and \cite{JTr} for certain classes of $2$-component links where the invariants are explicitly calculated.

A type of decomposition of a knot's front diagram known as a {\it normal ruling} provides a central connection between generating family and holomorphic techniques.    This notion arose independently in the work of the first author on the existence of augmentations \cite{F} and as a combinatorial invariant defined by Chekanov-Pushkar \cite{ChP} motivated by generating families.  Rulings provide a link between the Chekanov-Eliashberg invariant and generating families by combining the following two theorems:

\bigskip

\noindent {\bf Theorem 2.7} (\cite{F} \cite{FI} \cite{S}) {\it  The Chekanov-Eliashberg DGA of $\ell$ admits a ($\rho$-graded) augmentation if and only if the front diagram of $\ell$ admits a ($\rho$-graded) normal ruling.}

\bigskip

\noindent {\bf Theorem 2.4} (\cite{ChP}) {\it A Legendrian knot $\ell$ has a generating family if and only if the front diagram of $\ell$ has a graded normal ruling. }

\bigskip

Considering the two theorems in conjunction we see that it is possible to form integer graded linearized homology groups $H^\varepsilon_*(\ell)$ if and only if $\ell$ has a generating family.  This suggests a generating family interpretation for the groups $H^\varepsilon_*(\ell)$.  Theorem 5.3 shows that the set of possible (isomorphism classes) of generating family homology groups is contained in the set of linearized contact homology groups,  $\{ {\scr G}H_*(F) \} \subset \{ H^\varepsilon_*(\ell) \} $.    

\paragraph{Question.}  Does the reverse inclusion hold?  More precisely, given a graded augmentation $\varepsilon$ of the Chekanov-Eliashberg DGA of $\ell$ is it possible to find a linear at infinity generating family $F$ for $\ell$ so that $ {\scr G}H_*(F)  \cong  H^\varepsilon_*(\ell)  $?

\medskip

In \cite{SD} Sabloff established a duality theorem for the linearized contact homology groups $H^\varepsilon_*(\ell)$.

\noindent {\bf Theorem 6.1} (\cite{SD})  {\it If $\ell$ is a Legendrian knot and $\varepsilon : {\mathbf{ A}} \rightarrow {\bb Z}_2$ any graded augmentation, then we have

\[\begin{array} {rll} \dim_{{\bb Z}_
2}H^\varepsilon_k({\ell})&= \dim_{{\bb Z}_2}
H^\varepsilon_{-k}(\ell) & k\neq \pm 1 \\ \dim_{{\bb Z}_
2}H^\varepsilon_1(\ell)&= \dim_{{\bb Z}_2}
H^\varepsilon_{-1}(\ell) +1. &  \end{array}\]}

\bigskip

From Theorem 5.3 it follows that the dimensions of the generating family homology groups $ {\scr G}H_*(F)$ satisfy the same relations.  For the cases $k \neq \pm 1$, we show in section 6 how to interpret Sabloff duality from the generating family perspective as a version of Alexander duality.   

\smallskip

\paragraph{Overview of the rest of the paper:}  In section 2 we recall necessary notions from Legendrian knot theory.  The Chekanov-Eliashberg DGA is defined as well as augmentations, normal rulings, and generating families.  We review here the construction of a normal ruling from a generating family from \cite{ChP}.  This proves the forward implication of Theorem 2.4. Finally,  as a corollary of Theorem 2.4 and the results in \cite{R} we show that if a Legendrian knot has a generating  family then its Thurston-Bennequin number is maximal among Legendrian knots with the same underlying topological knot type. 

Section 3 gives a proof of the reverse implication in Theorem 2.4 via an explicit construction of a generating family for a Legendrian link admitting a normal ruling.  Although, this result is stated in \cite{ChP} as far as the authors know a proof has yet to appear in print.  

In section 4, a version of the `splash' construction is presented which we use for computations with the Chekanov-Eliashberg DGA.  We give a proof of Theorem 2.7 in which the similarity between the construction of a normal ruling from either a generating family or an augmentation is emphasized.  

The statement and proof of Theorem 5.3 occupies all of section 5.  The analogy between the proofs of Theorem 2.4 and 2.7 motivates the construction of an augmentation directly from a generating family used here.  To compute the generating family homology groups we use a fiber-wise version of the Morse complex based on bifurcation data in a $1$-parameter family of functions. 

 In conclusion, section 6 addresses the Sabloff duality.  In Theorem 6.2 we include a proof of Sabloff duality for the groups ${\scr G}H_*(F)$ purely from the generating family perspective.

\smallskip

\subsection{Acknowledgements}  The second author received support from NSF VIGRE Grant No. DMS-0135345 during portions of this
work.  The authors are grateful to AIM for hosting a workshop on ``Legendrian and transverse knots" where results from this work were presented.  Also, we thank the referee for many detailed comments which have helped to improve this paper.  Last but certainly not least, the second author thanks the first author for being so generous with his time and knowledge throughout the second author's years in grad school.

\section{Survey of known results}

\subsection{Legendrian curves and their projections}

A {\sl Legendrian curve} in ${\bb R}^3$ (with respect to the standard contact 
structure) is an immersed smooth curve $x=x(t),y=y(t),z=z(t)$ satisfying the 
equation $y\dot x-\dot z=0$ (that is, tangent to the distribution ${\scr C}
=\{y\, dx-dz=0\}$). The words Legendrian knots, Legendrian links, and 
Legendrian isotopy have obvious sense. There are two convenient projections of 
Legendrian curves. An $xz$ projection, or a {\sl front projection} of a 
generic Legendrian curve is a smooth curve with finitely many cusps but 
without vertical tangents. A front projection uniquely determines the 
Legendrian curve: the missing $y$ coordinate is reconstructed as the slope of 
the tangent line. A front projection of a Legendrian knot may have 
self-intersections but no self-tangencies (the latters would correspond to 
self-intersections of the Legendrian curve in space). An $xy$ projection of a 
Legendrian curve is smooth; it determines the curve up to a translation in the 
direction of the $z$ axis: the missing $z$ coordinate is reconstructed as 
$\int y\, dx$. An $xy$ projection of a closed Legendrian curve, in particular, 
of a Legendrian knot, is a self-intersecting smooth curve enclosing a region with zero signed  
area.

\subsection{Classical invariants}

There are two classical integer-valued Legendrian isotopy invariants of 
Legendrian knots (and links). The {\sl Thurston-Bennequin number} $TB(\ell)$ 
of a Legendrian knot $\ell$ is the linking number $\lk(\ell,\ell^+)$ where 
$\ell^+$ is a curve obtained from $\ell$ by a small shift in the direction of 
a normal within the distribution $\scr C$. (For knots, this number does not 
depend on the direction of the normal.) The {\sl rotation number} $R(\ell)$ of 
an {\sl oriented} Legendrian knot $\ell$ may be defined as the rotation number 
of its $xy$ projection. 

Both $TB(\ell)$ and $R(\ell)$ have simple description in terms of a front 
diagram. In this article, we will need this only for $R(\ell)$. Let $L$ be a 
front diagram of an oriented Legendrian knot $\ell$. Cusps break $L$ into 
non-self-intersecting parts, {\sl ``strands''}. Take one of the strands and 
attach to it some integer, $k$. Let us move from the chosen strand in the 
direction of the chosen orientation of $L$. Passing through a cusp to a new 
strand, we add 1 to our integer, if near the cusp the new strand is above the old one and 
subtract 1 otherwise. When we return to the initial strand, our integer 
becomes some $k'$; it is easy to see that $k-k'=2R(\ell)$. It is important 
that if $R(\ell)=0$, then every strand acquires a number (which we will call 
the {\sl index}), and of two strands forming a cusp, the index of the upper 
one is one more than the index of the lower one. These indices are defined up 
to simultaneous adding the same integer to all of them. If $R(\ell)\ne0$, then 
indices with similar properties are defined as residues modulo $2|R(\ell)|$.

\subsection{The Chekanov--Eliashberg DGA}

Consider a generic $xy$-diagram $\Gamma$ of a Legendrian knot $\ell$. Let $S$ 
be the set of all crossings of 
$\Gamma$, and let $\mathbf{A}=\mathbf{A}(\Gamma)$ be a free associative unital 
${\bb Z}_2$-algebra generated by $S$. At every crossing $s\in S$, the diagram 
forms four corners, of which we declare two positive and two negative: if you 
approach the crossing $s$ along the upper strand (that is, the strand with a 
bigger value of $z$), then the corner to the right of you is positive and the 
corner to the left of you is negative.

For every $n\ge0$ fix a convex planar domain $P_n$ bounded by a piecewise 
smooth curve with $n+1$ corners numerated counterclockwise as $v_0,v_1,\dots,
v_n$. For a crossing $s$, consider the set $I_n(s)$ of regular isotopy classes 
of orientation preserving immersions $f\colon P_n\to{\bb R}^2$ such that (1) 
$f(\partial P_n)\subset\Gamma$, (2) $f(v_0)=s$, (3) a neighborhood of $v_0$ 
covers a positive corner at $s$, and (4) for $i=1,\dots,n$, a neighborhood of 
$v_i$ covers a negative corner at $f(v_i)$. Put $I(s)=\cup_nI_n(s)$. The 
differential $d\colon{\bf A}\to{\bf A}$ is defined as a derivation such that 
$d(s)=\sum_{[f]\in I(s)}f(v_1)\dots f(v_n)$. It is proved in \cite{Ch} that 
$d^2=0$.

There exists a natural grading of $\mathbf A$ which assigns to each 
crossing $s$ a degree which is an integer, if $R(\ell)=0$, and a residue 
modulo $2|R(\ell)|$ otherwise. Here is the definition. 
At a crossing $s$ choose a path $\gamma_s$ which leaves $s$ along the overcrossing and follows the knot
diagram until it returns to the crossing along the understrand. Let $r(\gamma_s)$ be the number of counter-clockwise rotations made by the tanget vector to $\gamma_s$.  Here we should either assume that at $s$ strands meet at $90^\circ$ angles or round $r(\gamma_s)$ to the nearest odd multiple of $1/4$.
Degree of $s$ is given by
\[
|s| = 2 r(\gamma_s) - \frac{1}{2} \,\,\, (\mbox{mod} \,\, 2R(\ell) )
\]
and does not depend on the choice of $\gamma_s$.
With respect to this grading the 
differential $d$ has the degree $-1$. 

The following result is the main achievement of the Chekanov--Eliashberg 
theory (\cite{Ch},\cite{El}).

\begin{theorem}[Chekanov, Eliashberg]\label{chel}

The (graded) homology of $\mathbf{A}(\Gamma)$ is a Legendrian isotopy 
invariant of $\ell$.

\end{theorem}

Chekanov's paper contains the following, more precise statement.

\begin{theorem}[Chekanov]\label{stable}

The stable isomorphism type of $\mathbf{A}(\Gamma)$ is a Legendrian isotopy 
invariant of $\ell$.

\end{theorem}

Let us provide an explanation. A {\sl stabilization} of $\mathbf{A}(\Gamma)$ is 
obtained from $\mathbf{A}(\Gamma)$ by adding two generators $a,b$ with 
$\deg(a)=\deg(b)+1$ and extending the differential by the formulas $d(a)=b,\, 
d(b)=0$. The algebras $\mathbf{A}(\Gamma)$, $\mathbf{A}(\Gamma')$ are stably 
isomorphic if there are iterated stabilizations of each that are isomorphic. (It is clear 
that a stabilization does not affect the homology, so Theorem \ref{chel} follows 
from Theorem \ref{stable}.)

\begin{figure}[hbtp]
\centering
\includegraphics[width=2.6in]{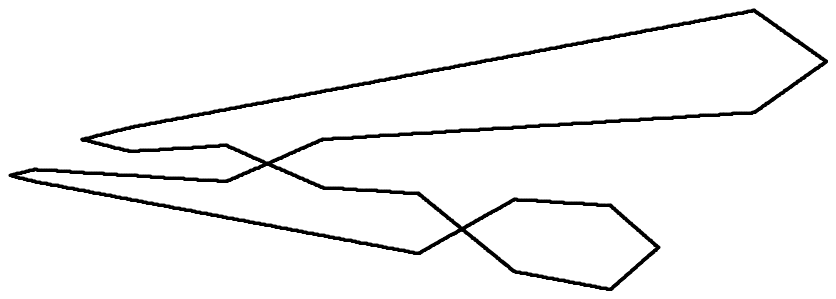}
\includegraphics[width=2.6in]{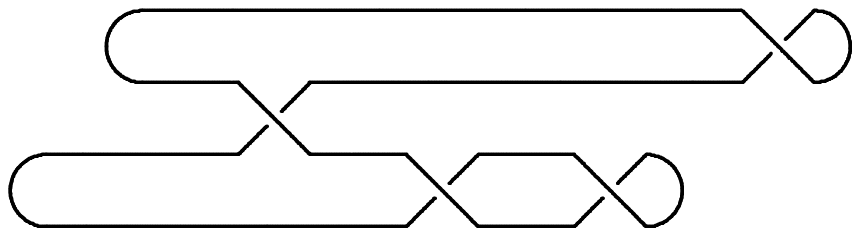}
\caption{Ng's construction}
\label{ngdiag}
\end{figure}

In conclusion, let us notice that a Chekanov--Eliashberg DGA of a Legendrian 
knot may be constructed from a front diagram. This was observed by L. Ng \cite{NgCLI} 
who showed that $xy$- and $xz$-diagrams of the same knot, however different 
they may seem, may be made to look almost the same.

A planar isotopy of the front diagram arranges the following.  Away from crossings and cusps the arcs of the front diagram are straight lines such that the slope of these lines decreases as we move from the top arc towards the bottom arc.  This will appear as several parallel horizontal lines on the $xy$-diagram.  At a left cusp the newly appearing arcs should have their slopes fit this criterion.  When we come to a crossing the slopes of the two adjacent arcs on the front diagram interchange causing an immediate crossing to appear on the $xy$-diagram and a crossing on the $xz$-diagram when the two lines eventually meet.  Before a right cusp the slopes of the two arcs that will meet are interchanged.  This produces an additional crossing on the $xy$-diagram which does not appear on the front diagram.

A diagram of this shape\footnote{The diagram in Figure 3 does not precisely follow Ng's construction.  Near crossings strands should change slope only once, but to accurately meet this criterion the diagram would need to become wider than this page.} 
is shown in Figure \ref{ngdiag}, left. The corresponding $xy$ diagram (Figure \ref{ngdiag}, right) looks, at 
least topologically, almost the same: crossings remain crossings, left cusps 
become roundings, right cusps become roundings preceded by additional 
crossings.

For this diagram, $\Gamma$, the algebra ${\bf A}(\Gamma)$ can be reconstructed 
from the initial front diagram, $L$, and we will denote it as ${\bf A}(L)$. 
The generators of ${\bf A}(L)$ correspond to crossings and right cusps of $L$. 
The degree of a generator corresponding to a crossing of strands $S,S'$ where 
$S$ has slope bigger than $S'$, is $\ind(S)-\ind(S')$ (ind stands for the 
index, see 2.2); the degree of a generator corresponding to a right cusp is 
$+1$.

\subsection{Rulings}

The notion of a ruling and its relationship on one hand with generating families 
and on the other with the Chekanov-Eliashberg DGA is a central motivation for this work. 
It was introduced in 2000 independently by Chekanov and Pushkar \cite{ChP} and the 
first author \cite{F}. 

Let $L$ be a front diagram of a Legendrian knot $\ell$. A {\sl ruling} of $L$ 
consists of 

\noindent (1) a correspondence between left and right cusps, and 

\noindent (2) for 
every pair of corresponding left and right cusps, $l$ and $r$, two disjoint 
(except $l$ and $r$) paths within the diagram and joining $l$ to $r$.  The paths should have strictly increasing 
$x$-coordinate, and paths joining different pairs of 
cusps can meet only at crossings.

Obviously, the paths of a ruling never pass through the cusps, except at the 
endpoints, and cover the whole diagram; this covering is one-fold, except at the 
crossings and the cusps, where it is two-fold. In particular, any crossing 
belongs to two paths which may exchange or not exchange the strands passing 
through the crossing. In the first case the crossing is called a {\sl switch}. 
A ruling is fully determined by the set of switches, so we can consider 
rulings as subsets of the set of crossings.  For each value of the $x$-coordinate not containing a crossing or cusp the division of the paths of the ruling
into pairs gives rise to a fixed point free involution of the strands of the front diagram.  The constructions of rulings in Theorems 2.4 and 
2.7 are given by describing this involution.

The notion of a ruling turns out to be useful only in the presence of the 
following ``{\sl normality condition}''. Assume that no two crossings of the 
diagram have the same $x$-coordinate. Let $s$ be a switch, let $p_u$ and 
$p_\ell$ be the upper and lower paths of 
the ruling passing through $s$, and let $q_u$ and $q_\ell$ be the other paths 
joining the same cusps as $p_u$ and $p_\ell$. Let $z, z_u$, and $z_\ell$ be 
the $z$-coordinates of $s$ and of the intersection points of $q_u$ and 
$q_\ell$ with the vertical line through $s$. We call the switch {\sl normal}, 
if $z_u>z>z_\ell,\ \mbox{or}\ z_\ell>z_u>z,\ \mbox{or}\ z>z_\ell>z_u.$ (In 
the remaining cases, $z_u>z_\ell>z,\ z_\ell>z>z_u,\ \mbox{and}\ 
z>z_u>z_\ell,$ the switch is {\sl abnormal}.)
A ruling is called normal if all the switches are normal.\smallskip

\noindent{\sl Example.} Let $L$ be the front diagram in Figure 2 (this is a 
``standard'' trefoil knot), with crossings $s_1,s_2,s_3$. Then there are four 
rulings, $\{s_1,s_2,s_3\},\{s_1\},$ $\{s_2\},\{s_3\}$, and only one of them, 
$\{s_2\}$, is not normal.

\begin{figure}[hbtp]
\centering
\includegraphics[width=1.8in]{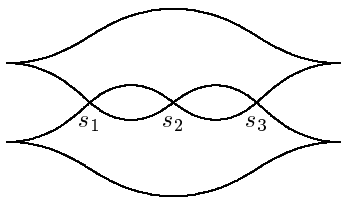}
\caption{Front diagram of a standard trefoil}
\label{trefoil}
\end{figure}

In the case when $R(L)=0$, a ruling is called {\sl graded}, if all the 
switches have degree 0. For example, all rulings of the diagram in Figure 
\ref{trefoil} are graded. In the general case, a ruling of $L$ is called 
$\rho$-graded, where $\rho$ is a divisor of $2R(L)$, if the degrees of all 
switches are divisible by $\rho$.  It is proved in \cite{S} that if $R(L)\ne0$, then 
there never exists a 2-graded normal  ruling.

\begin{figure}[hbtp]
\centering
\includegraphics[width=1.6in]{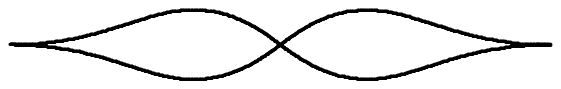}
\hskip.4in
\includegraphics[width=1.6in]{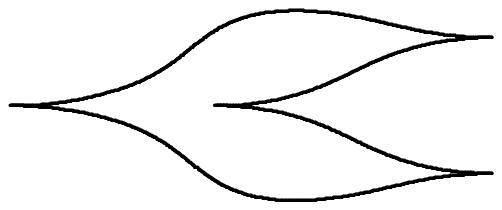}
\caption{Front diagrams with no rulings}
\label{noruling}
\end{figure}

Not every front diagram has a ruling (see Figure \ref{noruling}). But it is 
known \cite{ChP} that not only the existence of a normal ruling, but even the 
number of different normal rulings, as well as the number of different graded, 
or $\rho$-graded, normal rulings, is a Legendrian isotopy invariant. Moreover, a generic 
Legendrian isotopy between two front diagrams gives rise to a canonical bijection 
between the sets of normal rulings of these two diagrams. For further statements 
of this kind, see section 2.7.

There are at least three different theorems stating that a certain property of 
a Legendrian knot is equivalent to the existence of a (in some cases, graded) 
normal ruling. We will observe these results in the three subsections below. 
These theorems show, in particular, that some, visibly unrelated, properties of 
front diagrams are mutually equivalent. Below, we will establish a relation 
between two of them.

\subsection{Normal rulings and generating families} 

Let $F: \R^N \times \R \rightarrow \R$ be a smooth $1$-parameter family of functions, and let $f_t(x) = F(x,t), t\in{\bb R},\ x=(x_1,\dots,x_N)\in{\bb R}^N$.
Assume as well that $0$ is a regular value for the map
\[
(\frac{\partial F}{\partial x_1},\ldots, \frac{\partial F}{\partial x_N }): \R^N \times \R  \rightarrow  \R^N.
\]
Then, the fiber-wise critical set
\[
S_F = \{ (x,t) | (df_t)_x= 0 \}
\]
is a $1$-dimensional submanifold which becomes immersed in $\R^3$ as a Legendrian submanifold $\ell_F$ according to the mapping
\[
S_F \ni (x,t) \mapsto \left( t, \frac{\partial F}{\partial t} (x,t), F(x,t) \right).
\]
$F$ is called a {\it generating family} for $\ell_F$.  We require in addition that a generating family $F$ be {\it linear at infinity}.  That is, there exists a non-zero linear function $l:\R^N\rightarrow \R$  and $C >0$ such that $f_t(x)=l(x)$ when $|t|>C$ or $|x|>C$.  We will often take $l(x) = x_N$ or $l(x) = x_1$.  Clearly, any linear at infinity generating family can be given this form by precomposing with an appropriate linear transformation of $\R^N$.

Note that the front diagram $L$ of $\ell_F$ may be easily described as $$L=\{(t,z)\in{\bb R}^2\mid z\ 
\mbox{is a critical value of}\ f_t\}.$$
For a generic generating family $F$, for all but finitely many $t$ (which we will call exceptional),  $f_t$ is a Morse function with all critical values being 
different, and each exceptional value of $t$ corresponds to one of the 
following events: a generic birth or death of a pair of critical points of 
adjacent indices; a generic collision of two critical values.  Exceptional values of $t$ appear on the front diagram $L$ as cusps and crossings.

Which front diagrams possess generating families? One restriction is obvious.

\begin{proposition} If $L$ possesses a generating family, then $R(L_i)=0$ for each component $L_i$ of $L$. 
\end{proposition}

\paragraph{Proof} For any strand $S=\{(t,s(t)), t'\le t\le t''\}$ of $L_i$, 
there exists a continuous family $\{x_t\in{\bb R}^n, t'\le t\le t''\}$ where 
$x_t$ is a critical point of $f_t$ and $f_t(x_t)=s(t)$; obviously, the index 
of the critical point $x_t$ is the same for all $t$, and we denote this index 
by $\ind(S)$. At cusps of $L_i$, these indices behave as prescribed in 
Subsection 2.2; thus $R(L_i)=0$.

\begin{theorem}[Chekanov, Pushkar] A front diagram $L$ with $R(L)=0$ possesses 
a  (linear at infinity) generating family, if and only if it possesses a graded normal ruling.
\end{theorem}

The {\sl only if} part of this theorem is proved in \cite{ChP}. For the reader's convenience, we recall this construction of 
a normal ruling for a front diagram with a generating family of functions. Let $f\colon{\bb R}^N\to{\bb R}$ be a Morse function with all 
critical values different and with $f(x)=x_N$ for $|x|>C$. For a real $c$, set 
$X_c=\{x\in{\bb R}^N\mid f(x)\le c\}$. Call critical values $c_1>c_2$ of $f$ 
related, if, for a small $\varepsilon>0$,\[\begin{array} {rcl} \dim_{{\bb Z}_
2}H_\ast(X_{c_1+\varepsilon},X_{c_2-\varepsilon};{\bb Z_2})+1&=&\dim_{{\bb Z}_2}
H_\ast(X_{c_1-\varepsilon},X_{c_2+\varepsilon};{\bb Z_2})+1\\ =\dim_{{\bb Z}_2}
H_\ast(X_{c_1+\varepsilon},X_{c_2+\varepsilon};{\bb Z_2})&=&\dim_{{\bb Z}_2}
H_\ast(X_{c_1-\varepsilon},X_{c_2-\varepsilon};{\bb Z_2}).\end{array}\] It 
turns out that the pairs of related critical values are disjoint; moreover, 
the whole set of critical values falls into the union of disjoint related 
pairs; moreover, the involutions arising in the sets of critical values of the 
functions of a generating family compose a graded normal ruling of the front 
diagram. (In particular, related critical values have adjacent indices.)

A different way to describe this ruling arises from the following proposition 
(which is used in topology since J.H.C.Whitehead's works of the 1930's). 
Note that this proposition will be important for us in the subsequent parts 
of this article.

\begin{proposition} Let $V$ be a vector space over some field with a fixed 
ordered basis $e_1,\dots,e_m$, and let $d\colon V\to V$ be a linear 
transformation with the following two properties: $(1)$ $d$ is triangular, in 
the sense that $d(e_i)=\sum_{j>i}a_{ij}e_j$; $(2)$ $d$ is exact\begin{footnote} {There is a version of this proposition for the more general case $d^2=0$, but we do not need it.}\end{footnote} in the sense 
that $\Ker d=\Im d$. Then there exists a fixed point free involution 
$\tau\colon\{1,\dots,m\}\to\{1,\dots,m\}$ and a triangular basis change 
$e'_i=e_i+\sum_{j>i}b_{ij}e_j$ such that $d(e'_i)=e'_{\tau(i)}$, if 
$\tau(i)>i$, and $d(e'_i)=0$, if $\tau(i)<i$. Moreover, the involution $\tau$ 
with these properties is uniquely determined by $d$.\end{proposition}

For a non-exceptional value of $t$, consider the ${\bb Z}_2$ Morse complex
associated with the function $f_t$ and some Riemannian metric on ${\bb R}^N$ 
(compatible with $f_t$ in the sense of Morse-Smale). 
The total space $V$ of this complex has a natural basis corresponding to the 
critical values of $f_t$ ordered accordingly to the decreasing order in the 
set of critical values, and the differential $d$ is triangular (obviously) and 
exact (because the morse complex computes the homology of the pair $({\bb R}^N,\{f_t<-C\})$). 
According to Proposition 2.5, this provides an involution in the set of 
critical values, and it can be checked that it is the same involution as 
before (in particular, it does not depend on the Riemannian metric).  The involutions can be seen to piece together to give a graded normal ruling.  
In particular, the normality condition may be seen to hold as in the work of Barannikov \cite{B} (also see the proof of Proposition~\ref{norm}).

The {\sl if} part of Theorem 2.4 was proved by Pushkar, but, as far as we 
know, the proof has never been published. We provide a proof in Section 3. Our 
proof contains an explicit construction of a generating family for a front 
diagram with a graded normal ruling. 

\subsection{Normal rulings and augmentations} An augmentation of the 
Che\-kanov-Eliashberg DGA, $\mathbf{A}=\mathbf{A}(\Gamma)$, is a unital ring homomorphism 
$\varepsilon\colon \mathbf{A}(\Gamma)\to{\bb Z}_2$ such that $\varepsilon\circ d=0$. An 
augmentation is completely determined by its restriction to the set of 
generators of $\mathbf{A}$, that is, the set $S$ of crossings of the $xy$-diagram 
$\Gamma$ (the set of crossings and right cusps of the front diagram $L$ for 
$\mathbf{A}=\mathbf{A}(L)$). An augmentation is called graded ($\rho$-graded), if 
$\varepsilon(s)\ne0,\ s\in S$, implies $\deg s=0\ (\deg s\equiv0\bmod\rho)$. 

For an augmentation $\varepsilon$, set $\mathbf{A}^\varepsilon=\Ker\varepsilon/(\Ker
\varepsilon)^2$. This is a vector space with the basis $\{s^\varepsilon=
s+\varepsilon(s)\mid s\in S\}$. If $a\in\Ker\varepsilon$, then $da\in\Ker
\varepsilon$ (actually, $da\in\Ker\varepsilon$ for any $a\in \mathbf{A}$), and if 
$a\in(\Ker\varepsilon)^2$, then $da\in(\Ker\varepsilon)^2$ (if $a=bc,\ b,c\in
\Ker\varepsilon$, then $da=(db)c+b(dc)\in(\Ker\varepsilon)^2$). Hence, 
$d\colon \mathbf{A}\to \mathbf{A}$ gives rise to a homomorphism $d^\varepsilon\colon 
\mathbf{A}^\varepsilon\to \mathbf{A}^\varepsilon$, and $(d^\varepsilon)^2=0$. If the 
augmentation $\varepsilon$ is graded ($\rho$-graded), then $\mathbf{A}^\varepsilon$ is 
graded ($\rho$-graded), and $d^\varepsilon$ has the degree $-1$. 

$(\mathbf{A}^\varepsilon, d^\varepsilon)$ is refered to as the {\it linearized complex} of $\mathbf{A}$ with respect to the augmentation $\varepsilon$.  We denote the corresponding {\it linearized homology} as $H^\varepsilon_*(\ell)$.   $H^\varepsilon_*(\ell)$ is graded ($\rho$-graded), if so is 
$\varepsilon$.

\begin{theorem}[Chekanov \cite{Ch}] The set of all (all graded, all 
$\rho$-graded) homologies $\{ H^\varepsilon_*(\ell) \}$ corresponding to all (all graded, 
all $\rho$-graded) augmentations $\varepsilon$ of $\mathbf{A}$ is a Legendrian isotopy invariant.
\end{theorem}

The problem of existence of an augmentation is solved by the following result.

\begin{theorem}[Fuchs, Ishkhanov, Sabloff] The algebra $\mathbf{A}(\Gamma)$ possesses a 
(graded, $\rho$-graded) augmentation if and only if the corresponding front 
diagram possesses a (graded, $\rho$-graded) normal ruling.\end{theorem}

The {\sl if} part of this theorem was proved by the first author in \cite{F} (and 
this was one of the initial motivations for the notion of a normal ruling). The 
{\sl only if} part was proved in \cite{FI} and, independently, in \cite{S}. We will 
discuss the proofs in Section 4.  In particular we give a new approach to the forward
 implication based on Proposition 2.5.

\subsection{Normal rulings and estimates for the Thurston-Bennequin 
number}  Within any fixed topological knot type $K$ there exist Legendrian knots $\ell\in K$ with $TB(\ell)$ negative of arbitrarily large magnitude.  However, the set $\{TB(\ell) | \ell\in K\}$ is bounded from above.  For instance, there are estimates in terms of the two-variable knot polynomials (\cite{FT}, \cite{CG}, \cite{Ta}, \cite{NgSABE}):

$$TB(\ell) \leq -\deg_aF_K -1 \eqno(1)$$
$$TB(\ell) + |R(\ell)| \leq -\deg_aP_k -1$$

\noindent where $F_K, P_K\in {\bb Z}[a^{\pm1}, z^{\pm1}]$ denote the Kauffman and HOMFLY polynomials respectively (see \cite{R} for the particular conventions).  It was conjectured in \cite{F} that a $1$-graded normal ruling exists if and only if the estimate (1) is sharp.  This follows from a stronger relationship.

\begin{proposition}[\cite{R}]  

The coefficient of $a^0$ in $a^{TB(\ell)+1}F_K$ (resp. $a^{TB(\ell)+1}P_K$) is given by $\sum_{r} z^{j(r)}$ where the sum is over all 1-graded (resp. 2-graded) normal rulings $r$ of $\ell$ and $$j(r) =\#\{\mbox{switched\ crossings}\}- \#\{\mbox{right\ cusps}\}+1.$$

\end{proposition}

It was shown in \cite{ChP} that in general the $\rho$-graded {\it ruling polynomial }, $R_\ell^\rho(z) = \sum_{r} z^{j(r)}$ where the sum is over all $\rho$-graded rulings, is a Legendrian isotopy invariant.  Proposition 2.8 shows that in the cases $\rho = 1,2$, $R^\rho$ depends only on $TB(\ell)$ and the underlying topological knot type.  In contrast, $R^0$ can distinguish between knots with identical classical invariants \cite{ChP}.

The following is a simple corollary of Theorem 2.4 and Proposition 2.8.

\begin{corollary}

If $\ell$ possesses a generating family then $TB(\ell)$ is maximal within the underlying topological knot type of $\ell$.

\end{corollary}

\paragraph{Proof.}  If $\ell$ has a generating family then, according to Theorem 2.4, its front diagram has a normal ruling.  Therefore, by Proposition 2.8, the coefficient of $a^0$ in $a^{TB(\ell)+1}F_K$ is non-zero, and it follows that the estimate (1) is sharp.

\section{Construction of a generating family of functions for a front diagram 
with a ruling.}

Let $L$ be a front diagram in the $xz$ plane equipped with a graded 
normal ruling $R$.
The condition of being graded includes an assumption that for 
every strand $S$ of $L$, an integer $\ind(S)$ (called the index of $S$) is 
assigned such that of two strands meeting at a cusp, the upper one has the
index one more than the lower one.  At switches the two crossing strands have the same index.
We will  assume the index of strands to be chosen
so that $\ind(S) \geq 2$ for every $S$ (otherwise add the same constant to all of them). Choose 
a number $N$ so that $N-1$ exceeds all the indices.


We will construct a family of functions, $F\colon{\bb R}^N\times \bb R \to\bb R$, $f_t = F(\cdot, t)$,
 such that $L$ and $R$ are the diagram of critical values and the Morse 
complex ruling corresponding to $\{f_t\}$ in the sense explained in 2.5.  The index of strands will agree with the Morse index of critical points.  For values of $t$ away from switches, the $f_t$ will have a relatively simple form with critical points divided into pairs in accordance with the ruling $R$.  The necessary modifications at switches are presented below via description of the level surfaces.  All references to gradient trajectories and Morse complexes make use of the Euclidean metric.

The 
conditions imposed above on the indices guarantee that  functions $f_t$ will 
not have critical points of index $0$, $1$, $N-1$, and $N$.  It follows that the level sets of an $f_t$ immediately preceding and following a critical value will share the same $0$-th homology group.  In particular, since sufficiently negative level sets are simply hyperplanes it follows that all the level 
surfaces $f_t=a$ are connected.

Outside of a compact subset of its domain 
we will assign $F(x,t) = x_N$.  
Our family will have an additional property that for $t$ not equal to the $x$ 
coordinates of cusps and not belonging to small neighborhoods of the $x$ 
coordinates of switches of $R$, the Morse complex of $f_t$ has the simplest 
possible structure. Namely, the critical points of $f_t$ are arranged into 
pairs, according to the ruling, and in every such pair, the indices of 
critical points differ by one; we state that for every pair there will be 
precisely one gradient trajectory joining these two points, and there will be 
no other gradient trajectories joining the critical points of adjacent indices.  

For $t<-K$ (where $K$ is a positive number such that $L$ is contained in the 
domain $|t|<K$), we put $f_t(x_1,\dots,x_N)=x_N$. Moving $t$ to the right, we 
reach the leftmost cusp of $L$. At this moment, we create a pair of critical 
points of appropriate indices with appropriate critical values such that the 
point with the smaller value (= the smaller index) lies precisely above the 
other point (here and below, we regard the direction of the $x_N$ axis as 
vertical and upward).

\begin{figure}[hbtp]
\centering
\includegraphics[width=5.2in]{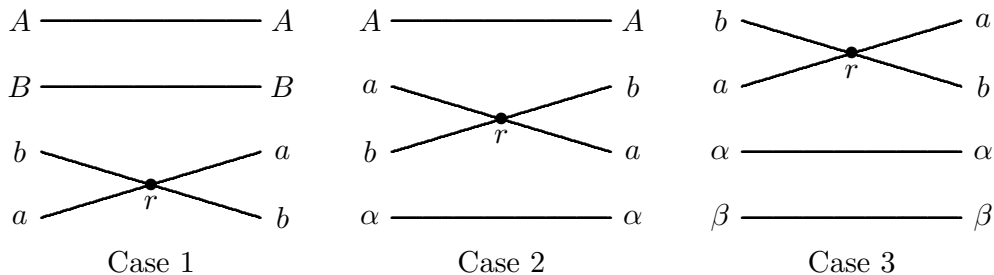}
\caption{Three cases for a switch of a ruling}
\label{normalswitch}
\end{figure}

 Moving further, we create, in a similar way, pairs of 
critical points at every left cusp in such a way that every new pair is 
located at a big distance (in terms of $x_1,\dots,x_{N-1}$ from the previous 
pairs.    
 When $t$ grows, but does not reach an $x$-coordinate of a switch from 
$R$, the critical values change according to the $z$-coordinates of points of 
$L$, and the pairs of critical points remain on lines parallel to the $x_N$ 
axis, and the gradient trajectories joining critical points of adjacent 
indices are all straight vertical lines. Nothing changes at non-switch crossings  
of $R$ (no crossings involve points from the same pair). If moving 
$t$ to the right we arrive at a right cusp, then the strands meeting at the 
cusp correspond to each other with respect to the ruling $R$; hence the 
corresponding critical points form a removable pair, and we remove them. It 
remains to explain what happens when we arrive at a crossing $r$ which is a switch for  
$R$. According to the definition given in 2.4, three cases are possible (see 
Figure \ref{normalswitch}.)

In Case 1, fragments of pairs of paths forming the ruling (directed from the 
left to the right) are $arb, AA$ and $bra, BB$; in Case 2 these fragments are 
$arb, AA$ and $bra, \alpha\alpha$ and in Case 3 $arb, \alpha\alpha$ and $bra, 
\beta\beta$. In all cases, $a$ and $b$ have the same index, say, $k$; then 
$A$ and $B$ have index $k+1$, and $\alpha$ and $\beta$ have index $k-1$. 

A family of functions corresponding to Case 1 is presented in Figure 5. The 
upper diagram shows the function $F_{t'}$ for some $t'<t$. The critical points 
$A$ and $a$ are joined by a gradient trajectory, and so are the critical 
points $B$ and $b$. The differential of the Morse complex acts as $A\mapsto 
a,\ B\mapsto b$. Then the two pairs of critical points are moved to each 
other.  At some instant a ``handle slide" trajectory appears when a gradient trajectory leaving $b$ ends up at $a$. 
This causes the Morse complex to change to $A\mapsto a,\ B\mapsto b+a$ as it appears in the next diagram.  (See Lemma 5.2 for a thorough discussion of handle slides and there effect on the Morse complex.)
Modulo a triangular 
transformation this is the same as before: the critical value at $b$ exceeds the critical value at $a$. 
Then these critical values are swapped (it is the crossing (the middle diagram 
corresponds to the function $F_t$ where $t$ is the $x$ coordinate of $r$), and 
the differential of the Morse complex (the same as before), $A+B\mapsto b,\ 
B\mapsto a+b$, becomes, after a triangular transformation, $A\mapsto b,\ 
B\mapsto a$ (now, the critical value at $a$ exceed that at $b$). 

\begin{figure}[hbtp]
\centering
\includegraphics[width=1.6in]{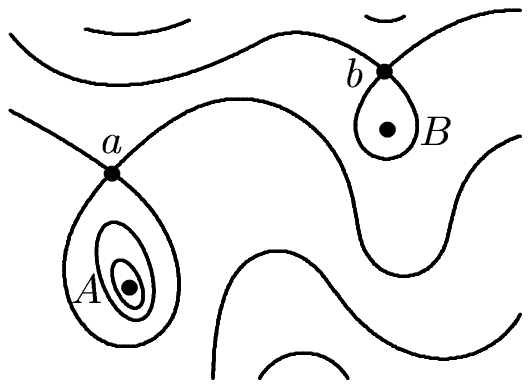}\hskip.4in\includegraphics[width=1.7in]{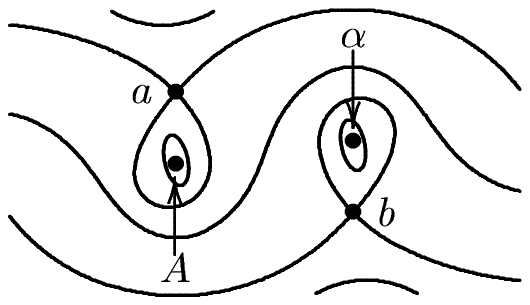}
\vskip.2in
\includegraphics[width=1.6in]{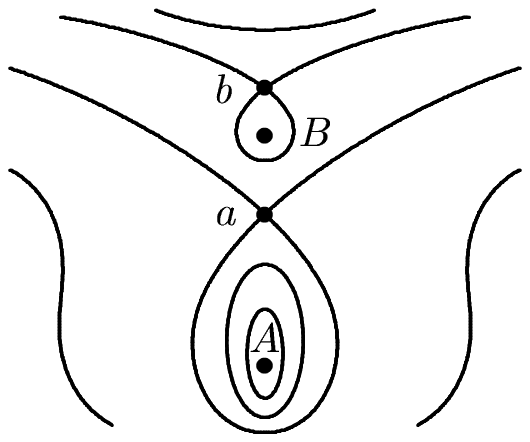}\hskip.4in\includegraphics[width=1.7in]{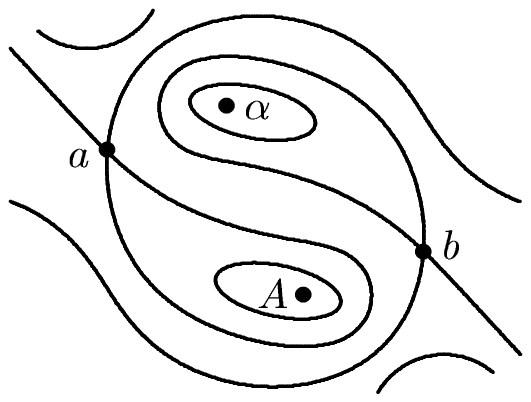}
\vskip.2in
\includegraphics[width=1.6in]{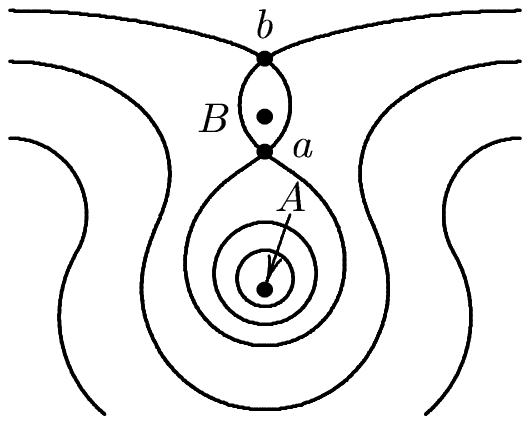}\hskip.4in\includegraphics[width=1.7in]{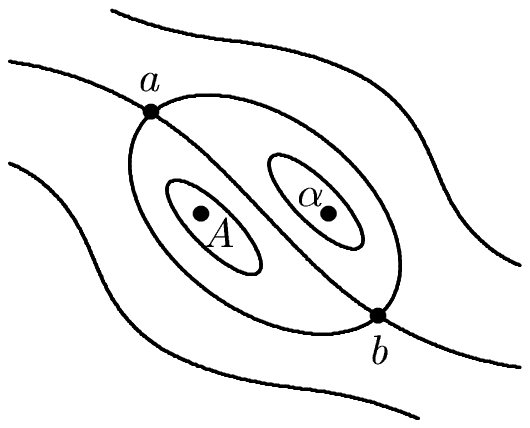}
\vskip.2in
\includegraphics[width=1.6in]{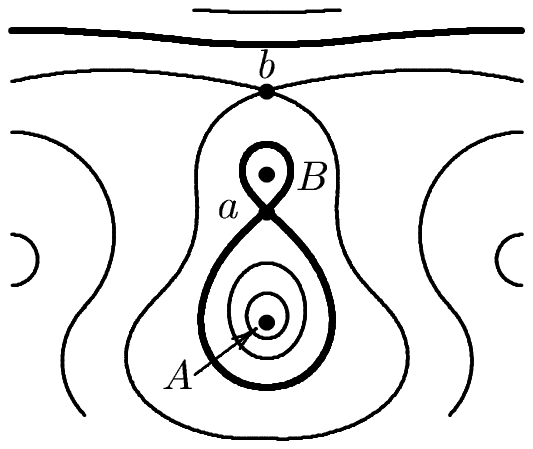}\hskip.4in\includegraphics[width=1.7in]{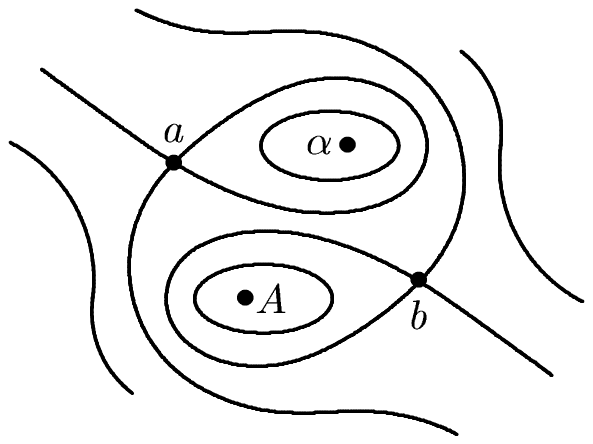}
\vskip.2in
\includegraphics[width=1.6in]{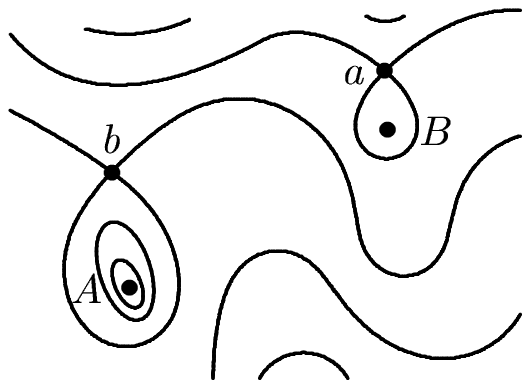}\hskip.4in\includegraphics[width=1.7in]{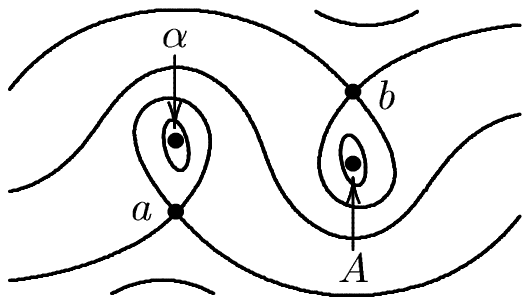}
\caption{Family of functions (Case 1 on the left, case 2 on the right)}
\label{generating}
\end{figure}

At last, the 
two new pairs of critical points, $A,b$ and $B,a$ are moved apart. (By the 
way, the transition from the second from the bottom diagram in Figure \ref{generating}, left
to the bottom one may seem not clear in Figure \ref{generating}, left; 
the reason is that our drawings are 2-dimensional, while the actual number 
of variables is at least 3; our assumptions on Morse indices have guaranteed the level surface of the point $a$ shown as a 
thick curve at the second from the bottom diagram is, actually, connected, and 
we can move the critical points $a$ and $B$ along this surface in an arbitrary 
way.  Handle slide trajectories connecting $a$ to $b$ as well as $A$ to $B$ will occur during this transition.) The bottom diagram corresponds to the function $f_{t''}$ with a $t''>t$; 
here, again, all gradient trajectories between critical points of adjacent 
indices are vertical lines. Remark, that our functions may have many other 
(pairs of) critical points, but all of them stay frozen in our deformation.

In Case 2, the family of functions is shown in Figure \ref{generating}, right. The critical points 
$A,a,b,\alpha$ have indices $k+1,k,k,k-1$. In the first diagram (corresponding 
to a $t'<t$), the Morse differential acts as $A\mapsto a,\ b\mapsto\alpha$. 
Then we move the two pairs of critical points to each other, and the 
differential becomes $A\mapsto a+b\mapsto 0, b\mapsto\alpha$ (a handle slide of $b$ over $a$ has occured) which is 
triangular equivalent to the previous differential while the value at $a$ 
exceeds the value at $b$. When these values are swapped (this happens when we 
pass from the second diagram to the forth one), then the differential is unchanged but becomes 
triangular equivalent to $A\mapsto b, a\mapsto\alpha$.  It remains only to move 
the pairs $A,b$ and $a,\alpha$ of critical points apart (see the bottom 
diagram which corresponds to the function $F_{t''}, t''>t$). During this last step a handle slide of $b$ over $a$ results in the Morse complex returning to a simple form  $A\mapsto b, a\mapsto\alpha$.

Case 3 is symmetric to Case 1 and is also illustrated by Figure \ref{generating}, left
(it is sufficient to reflect all the diagrams in horizontal 
lines and replace $A$ and $B$ by $\alpha$ and $\beta$).

\section{Relations between rulings, augmentations, and the linearized complex}

\subsection{Splash construction} Proofs of many results of the Legendrian 
knot theory, including Theorem 2.7, depend on (various versions of) the 
``splash construction'' which first appeared in \cite{F}. 
The goal of this construction 
is to modify an $xy$-diagram of a Legendrian knot in such a way that the 
differential of the Chekanov-Eliashberg DGA is described by explicit formulas; 
this is acheived at the expense of increasing the number of crossings. The many new generators may be organized into a sequence of matrices.  This perspective is useful later for  the construction of an augmentation from a generating family in Section 5.

\begin{figure}[hbtp]
\centering
\includegraphics[width=3.5in]{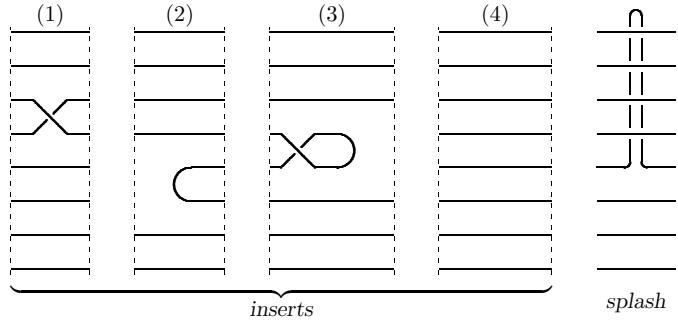}
\caption{Inserts and splashes.}
\label{insertsplash}
\end{figure}

We 
describe here a version of the splash construction.  Begin by applying Ng's resolution procedure (Section 2.3) to a  Legendrian knot so that the $xz$ and $xy$ diagram  have similar appearance.  The $xy$-diagram is cut by vertical lines 
into ``laminated zones'' separated by ``crossings'' and ``cusps''. Each 
laminated zone consists of several horizontal strands stretched between the 
vertical boundaries of the zone. The laminated zones are separated by inserts 
of four types (``crossings,'' ``left cusps,''  ``right cusps,'' or ``parallel lines'') shown in 
Figure \ref{insertsplash}, left.

\begin{figure}[hbtp]
\centering
\includegraphics[width=5.2in]{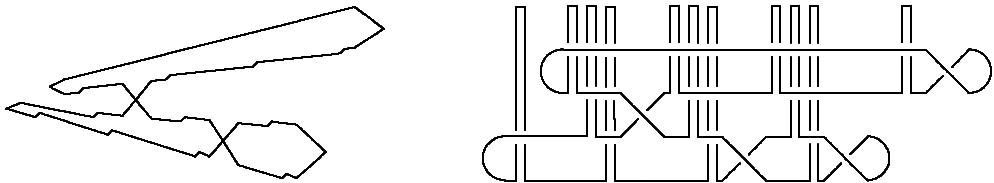}
\caption{Splash construction.}
\label{splashdiag}
\end{figure}

In the ``laminated zones'' we make ``splashes'' shown in Figure 
\ref{insertsplash}, right. We add a splash to each strand in the zone starting 
from the second top strand and ending with the bottom strand. Figure 
\ref{splashdiag}, right, shows the diagram of Figure \ref{ngdiag}, right, 
modified by splashes. As to the $xz$-diagram, the splashes leave it almost 
unchanged: low steep steps appear on its strands (Figure \ref{splashdiag}, 
left).

\begin{figure}[hbtp]
\centering
\includegraphics{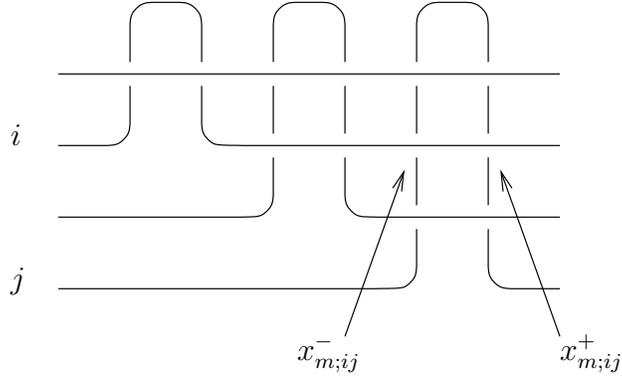}
\caption{A labeling of the strands produces a labeling of the crossings in laminated zones.}
\label{CL}
\end{figure}

To describe the Chekanov-Eliashberg DGA of the splashed diagram, we need 
notations for the crossings in the latter. Numerate the laminated zones from 
the left to the right by the numbers $1,2,\dots,N$, and the inserts, also from 
the left to the right, by the numbers $1,2,\dots,N+1$; thus, the insert to the 
left of the laminated zone number $m$ has also number $m$. The strands within 
a laminated zone we label, from the top to the bottom, by elements of an 
ordered set, usually (but not always, see below), by the numbers 
$1,2,\dots,n$. The crossings of the $m$-th splash between the $j$-th strand and the 
$i$-th strand $(i<j)$ we denote as $x^\pm_{m;ij}$, with $x^-_{m;ij}$ to the 
left of $x^+_{m;ij}$  (see Figure \ref{CL}). The crossing within $m$-th insert we denote as $y_m$, if 
this crossing arises from a crossing of the front diagram, and as $z_m$, 
if this crossing arises from a right cusp.

Adding splashes in the laminated zones has a simplifying effect on $d$.  The polygons that we need to count become trapped between two adjacent laminated zones so that $d$ of a generator coming from the $m$-th laminated zone is contained in the sub-algebra generated by generators from the $m$-th and $(m-1)$-st laminated zones and $m$-th insert.  

\begin{figure}[hbtp]
\centering
\includegraphics{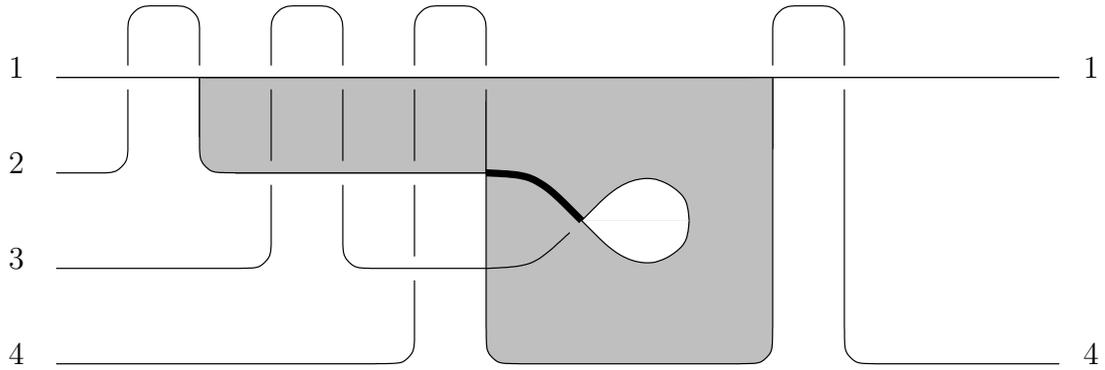}
\caption{The pictured disk accounts for an $x^+_{m-1;12} z_m x^+_{m-1;24}$ term in $d x^{-}_{m;14}$.}
\label{SplashDisk}
\end{figure}

\paragraph{ {\bf Labeling convention.}} The differential $dx^-_{m;ij}$ depends only on the $m$-th insert and consists of polynomials in the $x^\pm_{m;ij}$, $x^+_{m-1;ij}$ and $z_m$ or $y_m$.  Unfortunately, if the $m$-th insert contains a cusp then the convention of using consecutive integers,  $1,2, \dots, n$, to label strands 
can result in the same strand having distinct labels in the $m$-th and $m-1$-st laminated zones. This unnecessarily complicates the formulas for  $dx^-_{m;ij}$.
Instead, when the $m$-th insert contains a cusp we will use an alternate labeling of strands for the $m$-th and $m-1$-st when  discussing 
$d x^-_{m;ij}$ as follows:

If  the $m$-th insert contains a left (resp. right) cusp, then we numerate the strands in the $m$-th (resp. $m-1$-st) zone by the numbers $1,2,\dots,n$, and in the $(m-1)$-st (resp. $m$-th) zone by the same numbers except the numbers $k$ and $k+1$ corresponding to the strands meeting at the cusp are omitted. (See Figure \ref{SplashDisk}.)

\medskip

We warn the reader that this convention can cause the same crossing to be described using two different labels depending on the context.  
For instance, when discussing $d x^-_{m-1;ij}$ we use a (possibly different) labeling of the strands depending on what the $(m-1)$-st insert is. In the formula for $d x^+_{m;ij}$ we always use the standard labeling of $1, \ldots, n$ increasing from top to bottom.

\begin{proposition}

For all $m$, 
\[
dx^+_{m;ij}=\sum\limits_{i<s<j}x^+_{m;is}x^+_{m;sj}.
\]

{\rm (i)} If the $m$-th insert contains a crossing $y_m$ between $k$-th and 
$(k+1)$-st strands, then
\[\begin{array} {rcl} 
dy_m & =& x^+_{m-1;k,k+1},\quad \mbox{and}\\ dx^-_{m;ij}&=&x^+_{m;ij}+\sum\limits_{i<s<j}\left[
x^+_{m;is}x^-_{m;sj}+x^-_{m;is}\widetilde x^+_{m-1;sj}\right]+\widetilde 
x^+_{m-1;i,j},\end{array}\eqno(1)\]where $\widetilde x^+_{m-1;k,k+1}=0$, and, for 
$u<v, (u,v)\ne(k,k+1)$, \[\widetilde x^+_{m-1;uv}=x^+_{m-1;\bar u\bar v}+
\left\{\begin{array} {l} y_mx^+_{m-1;uv},\ \mbox{if}\ u=k+1,\\ 
x^+_{m-1;uv}y_m,\ \mbox{if}\ v=k,\end{array}\right.\]where, in turn, 
$\overline w=w$, if $w\ne k,k+1$, $\overline k=k+1,\overline{k+1}=k$.

{\rm(ii)} If the $m$-th insert contains a left cusp, then (using the labeling convention) the differentials $dx^-_{m;ij}$ are expressed by the same formulas (1) with a simpler formula for $\widetilde x^+_{m-1;uv}$:\[\widetilde 
x^+_{m-1;uv}=\left\{\begin{array} {ll} x^+_{m-1;uv},&\mbox{if}\ u\ne k,k+1,\ v\ne k,k+1,\\ 1,&\mbox{if}\ u=k,\ v=k+1\\ 0\strut &\mbox{otherwise}.\end{array}\right.\]

{\rm (iii)} If the $m$-th insert contains a right cusp with crossing $z_m$, then, (using the labeling convention) for $i\ne k,k+1,\ j\ne k,k+1$, 
\[\begin{array} {rcl}  
dz_m& =&  1+x^+_{m-1;k,k+1}, \mbox{and}\\ dx^-_{m;ij}&=&x^+_{m;ij}+\sum\limits_{i<s<j\atop s\ne k,k+1}\left[
x^+_{m;is}x^-_{m;sj}+x^-_{m;is}\widetilde x^+_{m-1;sj}\right]+\widetilde 
x^+_{m-1;i,j},\end{array}\]
where \[\begin{array}  {l} \widetilde x^+_{m-1;uv}=x^+_{m-1;uv}\end{array}\]
unless $u<k, k+1<v$ in which case

\[\begin{array} {lc} \widetilde x^+_{m-1;uv} =  x^+_{m-1;uv} + x^+_{m-1;u,k+1}x^+_{m-1;k,v} \  \\  \ + x^+_{m-1;uk}z_mx^+_{m-1;k,v}+ x^+_{m-1;u,k+1}z_mx^+_{m-1;k+1,v} + x^+_{m-1;uk}z_m^2x^+_{m-1;k+1,v} \end{array}\]

 {\rm (iv)}  If the $m$-th insert is simply $n$ parallel lines (does not correspond to any crossing or cusp on the front diagram), then

\[\begin{array} {rcl} dx^-_{m;ij}= & x^+_{m;ij} + \sum\limits_{i<s<j} x^+_{m;is}x^-_{m;sj} & + x^+_{m-1;ij} + \sum\limits_{i<s<j}x^-_{m;is}x^+_{m-1;sj} \end{array} \] 

\end{proposition}

Remark that in all cases the first formula has a short matrix form: if $X^+_m$ is the (strictly upper triangular) matrix with the entries $x^+_{m;ij}$, then $dX^+_m=(X^+_m)^2$. The other formulas in Proposition 4.1 also have matrix presentations; 
$$dX^-_m=X^+_m(I+X^-_m)+(I+X^-_m)\widetilde X^+_{m-1}$$
where the matrix $\widetilde X^+_{m-1}$ depends on the location and type of the singularity within the $m$-th insert.  For instance, if the insert is a crossing between the $k$ and $k+1$-st strands then $\widetilde X^+_{m-1} =A\widehat X^+_{m-1}A^{-1}$ where $\widehat X^+_{m-1}$ is the matrix with the entries $x^+_{m-1;ij}$ except for the entry $x^+_{m-1;k,k+1}$ which is removed, and $A$ being the block diagonal matrix with the block in the $k$-th and $(k+1)$-st rows and columns being $\displaystyle{\left[ \begin{array} {cc} 0&1\\ 1& y_m \end{array}\right]}$ and all other diagonal blocks being just $[1]$.

\subsection{From an augmentation to a normal ruling: a new approach} Let $L$ be a Ng front diagram of a Legendrian knot, and let $\Gamma$ be the $xy$-diagram obtained from $L$ by the splash construction (see Figures 1 and 8). We present here an algorithm which assigns to a (graded) augmentation of the Chekanov-Eliashberg DGA $\mathbf{A}(\Gamma)$ a (graded) normal ruling of $L$. This provides a new approach to the result of Fuchs, Ishkhanov and Sabloff (Theorem 2.7) and reveals its connection with Proposition 2.5 ( and hence with generating families).

$\mathbf{A}(\Gamma)$ is described in Proposition 4.1; below we use the notations of this Proposition. Let $\varepsilon: \mathbf{A}(\Gamma)\rightarrow \bb{Z}_2$ be an augmentation. Put $e^\pm_{m;ij}=\varepsilon(x^\pm_{m;ij})$, and let $E^\pm_m$ be the matrix with the entries $e^\pm_{m;ij}$ (that is, $E^\pm=\varepsilon(X^\pm_m)$). The matrix $E^\pm_m$ is triangular (that is, $e^\pm_{m;ij}=0$ for $j\le i$). Since $\varepsilon\circ d=0$, the relation $dX^+_m=(X^+_m)^2$ (see Section 4.1) shows that $(E^+_m)^2=0$.  The matrices $E^+_m$ will now play a role similar to that of the Morse complexes of the $f_t$ in Section 2.5.  

Let ${\scr C}_m$ be the vector space with the basis $\{\sigma_i\}$ labelled with the numbers of strands in the $m$-th laminated zone. Put $\partial(\sigma_i)=\sum e^+_{m;ij}\sigma_j$. Since $(E^+_m)^2=0$, $({\scr C}_m,\partial)$ is a complex (with a triangular differential, as in Proposition 2.5). Moreover, if the rotation number of the Legendrian knot is 0 and the augmentation $\varepsilon$ is graded, then ${\scr C}_m$ is graded and $\deg\partial=-1$.

\begin{lemma}

The complex ${\scr C}_m$ is acyclic.

\end{lemma}

\paragraph{Proof} It is sufficient to prove that the homologies of the complexes ${\scr C}_m$ and ${\scr C}_{m-1}$ are the same. Since ${\scr C}_m=0$ for $m<0$, this implies the acyclicity.

If the $m$-th insert contains a crossing $y_m$, then the formula $$dX^-_m=X^+_m(I+X^-_m)+(I+X^-_m)A\widehat X^+_{m-1}A^{-1}$$ (see Section 4.1) implies $$E^+_m=(I+E^-_m)
\varepsilon(A)\varepsilon(\widehat X^+_{m-1})\varepsilon(A)^{-1}(I+E^-_m)^{-1}.\eqno(2)$$But since 
$dy_m=x^+_{m-1;k,k+1}$ (see Part (i) of Proposition 4.1), $\varepsilon(x^+_{m-1;k,k+1})=0$, and hence $\varepsilon(\widehat X^+_{m-1})=E^+_{m-1}$ (the matrices $X^+_{m-1},\widehat X^+_{m-1}$ are the same with the exception of the entry $x^+_{m;k,k+1}$), and the formula shows that the matrices 
$E^+_m,E^+_{m-1}$ are conjugated. Hence the complexes ${\scr C}_m$ and ${\scr C}_{m-1}$ have the same homology.

The cases when the $m$-th insert contains the crossing $z_m$, or does not contain any crossing at all, are easier.

Let the $m$-th insert contain no crossings. The complex ${\scr C}_{m-1}$ has the basis $\{\sigma_i\mid i\ne k,k+1\}$. Let $(\widetilde{\scr C}_{m-1},\widetilde\partial)$ be $({\scr C}_{m-1},\partial)$ with basis elements $\sigma_k,\sigma_{k+1}$ added and $\widetilde\partial$ the same as $\partial$ with, additionally, $\widehat\partial(\sigma_k)=\sigma_{k+1}$. Then the formula (2) becomes$$E^+_m=(I+E^-_m)\widetilde E^+_{m-1}(I+E^-_m)^{-1}$$where $\widetilde E^+_{m-1}$ is the matrix of $\widetilde\partial$. Hence, the complex ${\scr C}_m$ has the same homology as the complex $\widetilde{\scr C}_{m-1}$ which, in turn, has the same homology as ${\scr C}_{m-1}$.

Let the $m$-th insert contain the crossing $z_m$. Then the complex ${\scr C}_{m-1}$ has two generators missing in ${\scr C}_m\colon\sigma_k$ and $\sigma_{k+1}$. Since $dz_m=1+x^+_{m-1;k,k+1}$ (see Part (iii) of Proposition 4.1), in ${\scr C}_{m-1}$, $\partial\sigma_k=\sigma_{k+1}+$ a linear combination of $\sigma_j$ with $j>k+1$. Consider the complex $\widetilde{\scr C}_{m-1}$ as the quotient of the complex ${\scr C}_{m-1}$ by the (acyclic) two dimensional subcomplex generated by \[\left\{ \begin{array} {lll} \sigma_k\ \mbox{and}\ \partial\sigma_k,&\mbox{if}&\varepsilon(z_m)=0,\\ \sigma_k+\sigma_{k+1}\ \mbox{and}\ \partial\sigma_k + \partial\sigma_{k+1},&\mbox{if}&\varepsilon(z_m)=1.\end{array}\right.\] We can assume that the complex $\widetilde{\scr C}_{m-1}$ has the same basis as ${\scr C}_m$. The last formula in Proposition 4.1 shows that $$E^+_m=(I+E^-_m)\widetilde E^+_{m-1}(I+E^-_m)^{-1}$$ and the matrix $\widetilde E^+_{m-1}$, whether  
$\varepsilon(z_m)=0$ or $1$, is the matrix of the differential of $\widetilde{\scr C}_{m-1}$ with respect to the basis $[\sigma_j]$, $j\neq k,k+1$. Hence, the complex ${\scr C}_m$ has the same homology as the complex $\widetilde{\scr C}_{m-1}$ which, in turn, has the same homology as ${\scr C}_{m-1}$.

This completes the proof of Lemma.\smallskip

According to Proposition 2.5, there arises a fixed point free involution $\tau_m$ in the set of generators of the complex ${\scr C}_m$, that is, in the set of strands in the $m$-th laminated zone of $L$, such that, after a triangular basis transformation $\sigma'_i=\sigma_i+\sum_{j>i}a_{ij}\sigma_j$, the differential $\partial\colon{\scr C}_m\to{\scr C}_m$ acts as \[\partial(\sigma'_i)=\left\{ \begin{array} {lll} \sigma'_{\tau_m(i)},&\mbox{if}&i<\tau_m(i),\\ 0,&\mbox{if}&i>\tau_m(i).\end{array}\right.\]

\begin{proposition}
\label{norm}
The involutions $\tau_m$ form a normal ruling of $L$, graded, if the augmentation $\varepsilon$ is graded.

\end{proposition}

This completes our algorithm: it provides a (graded) ruling from a (graded) augmentation.

\paragraph{Proof of Proposition 4.5} Let $m$-th insert contain a crossing $y_m$. After a triangular changes of basis, the complex ${\scr C}_{m-1}$ and ${\scr C}_m$ have differentials, respectively, \[
\partial\sigma_i=\left\{\begin{array} {ll} \sigma_{\tau_{m-1}(i)},&\mbox{if}\ \tau_{m-1}(i)>i,\\ 0,&\mbox{if}\ \tau_{m-1}(i)<i,\end{array}\right. \hskip.3in \partial\sigma_i=\left\{\begin{array} {ll} \sigma_{\tau_m(i)},&\mbox{if}\ \tau_m(i)>i,\\ 0,&\mbox{if}\ \tau_m(i)<i.\end{array}\right. \] The transition from ${\scr C}_{m-1}$ to ${\scr C}_m$ (up to an additional triangular transformation by the matrix $I+E^-_m$) is performed by the matrix $\varepsilon (A)$, that is,\[\sigma_i\mapsto\left\{\begin{array} {ll} \sigma_i,&\mbox{if}\ i\ne k,k+1,\\ \alpha\sigma_k+\sigma_{k+1},&\mbox{if}\ i=k,\\ \sigma_k,&\mbox{if}\ i=k+1,\end{array}\right.\]where $\alpha=\varepsilon(y_m)$. There are 6 cases.\smallskip

{\it Case} 1. $\tau_{m-1}(k+1)<\tau_{m-1}(k)<k$. In ${\scr C}_{m-1}$,$$\partial\sigma_{\tau_{m-1}(k+1)}=\sigma_{k+1},\, \partial\sigma_{\tau_{m-1}(k)}=\sigma_k;$$then in ${\scr C}_m$,$$\partial\sigma_{\tau_{m-1}(k+1)}=\sigma_k,\, \partial\sigma_{\tau_{m-1}(k)}=\alpha\sigma_k+\sigma_{k+1},$$ $$\partial(\sigma_{\tau_{m-1}(k+1)}+\sigma_{\tau_{m-1}(k)})=(1+\alpha)\sigma_k+\sigma_{k+1}.$$Hence,  if $\alpha=0$, then $\tau_m(k)=\tau_{m-1}(k+1),\, \tau_m(k+1)=\tau_{m-1}(k)$ (no switch), and if $\alpha=1$, then $\tau_m=\tau_{m-1}$ (a switch).\smallskip

{\it Case} 2. $\tau_{m-1}(k)<\tau_{m-1}(k+1)<k$. In ${\scr C}_{m-1}$,$$\partial\sigma_{\tau_{m-1}(k)}=\sigma_k,\, \partial\sigma_{\tau_{m-1}(k+1}=\sigma_{k+1};$$then in ${\scr C}_m$,$$\partial\sigma_{\tau_{m-1}(k)}=\alpha\sigma_k+\sigma_{k+1},\, \partial\sigma_{\tau_{m-1}(k+1)}=\sigma_k,$$ $$\partial(\sigma_{\tau_{m-1}(k)}+\sigma_{\tau_{m-1}(k+1)})=(1+\alpha)\sigma_k+\sigma_{k+1}.$$Hence, whether $\alpha=0$ or 1, $\tau_m(k)=\tau_{m-1}(k+1),\, \tau_m(k+1)=\tau_{m-1}(k)$ (no switch).\smallskip

{\it Case} 3. $\tau_{m-1}(k)<k,\, k+1<\tau_{m-1}(k+1)$. In ${\scr C}_{m-1}$,$$\partial\sigma_{\tau_{m-1}(k)}=\sigma_k,\, \partial\sigma_{k+1}=\sigma_{\tau_{m-1}(k+1)};$$then in ${\scr C}_m$,$$\partial\sigma_{\tau_{m-1}(k)}=\alpha\sigma_k+\sigma_{k+1},\, \partial(\alpha\sigma_k+\sigma_{k+1})=0,\, \partial\sigma_k=\sigma_{\tau_{m-1}(k+1)}.$$If $\alpha=0$, then $\tau_m(k)=\tau_{m-1}(k+1),\, \tau_m(k+1)=\tau_{m-1}(k)$ (no switch). If $\alpha=1$, then$$\partial\sigma_{\tau_{m-1}(k)}=\sigma_k+\sigma_{k+1},\, \partial\sigma_{k+1}=\partial\sigma_k=\sigma_{\tau_{m-1}(k+1)},$$hence $\tau_m=\tau_{m-1}$ (a switch).\smallskip

{\it Case} 4. $\tau_{m-1}(k+1)<k,\, k+1<\tau_{m-1}(k)$. In ${\scr C}_{m-1}$,$$\partial\sigma_{\tau_{m-1}(k+1)}=\sigma_{k+1},\, \partial\sigma_k=\sigma_{\tau_{m-1}(k)};$$then in ${\scr C}_m$,$$\partial\sigma_{\tau_{m-1}(k+1)}=\sigma_k,\, \partial(\alpha\sigma_k+\sigma_{k+1})=\sigma_{\tau_{m-1}(k)}.$$Hence,$$\partial\sigma_k=0,\, \partial\sigma_{k+1}=\sigma_{\tau_{m-1}(k)}.$$Thus, whether $\alpha=0$ or 1, $\tau_m(k)=\tau_{m-1}(k+1),\, \tau_m(k+1)=\tau_{m-1}(k)$ (no switch).\smallskip

{\it Case} 5. $k+1<\tau_{m-1}(k+1)<\tau_{m-1}(k)$. In ${\scr C}_{m-1}$,$$\partial\sigma_k=\sigma_{\tau_{m-1}(k)},\, \partial\sigma_{k+1}=\sigma_{\tau_{m-1}(k+1)};$$then in ${\scr C}_m$,$$\partial(\alpha\sigma_k+\sigma_{k+1})=\sigma_{\tau_{m-1}(k)},\, \partial\sigma_k=\sigma_{\tau_{m-1}(k+1)}.$$If $\alpha=0$, then $\tau_m(k)=\tau_{m-1}(k+1),\, \tau_m(k+1)=\tau_{m-1}(k)$ (no switch). If $\alpha=1$, then $\partial(\sigma_k+\sigma_{k+1})=\sigma_{\tau_{m-1}(k)}$, and hence $\tau_m=\tau_{m-1}$ (a switch).\smallskip

{\it Case} 6. $k+1<\tau_{m-1}(k)<\tau_{m-1}(k+1)$. In ${\scr C}_{m-1}$,$$\partial\sigma_k=\sigma_{\tau_{m-1}(k)},\, \partial\sigma_{k+1}=\sigma_{\tau_{m-1}(k+1)};$$then in ${\scr C}_m$,$$\partial(\alpha\sigma_k+\sigma_{k+1})=\sigma_{\tau_{m-1}(k)},\, \partial\sigma_k=\sigma_{\tau_{m-1}(k+1)}.$$If $\alpha=0$, then $\partial\sigma_{k+1}=\sigma_{\tau_{m-1}(k)},\, \partial\sigma_k=\sigma_{\tau_{m-1}(k+1)}$; if $\alpha=1$, then $\partial(\sigma_{k+1})=\sigma_{\tau_{m-1}(k)}+\sigma_{\tau_{m-1}(k+1)}.$ Hence, whether  $\alpha=0$ or 1, $\tau_m(k)=\tau_{m-1}(k+1),\, \tau_m(k+1)=\tau_{m-1}(k)$ (no switch).\smallskip

In addition to this, if the $m$-th insert contains no crossings, then in ${\scr C}_m$, $\partial\sigma_k=\sigma_{k+1}$, and if the $m$-th insert contains a crossing $z_m$, then in ${\scr C}_{m-1}$, $\partial\sigma_k=\sigma_{k+1}+$ generators with the numbers $>k+1$ (see the proof of Lemma 4.2). Thus, in the first case $\tau_m(k)=k+1$, and in the second case $\tau_{m-1}(k)=k+1$.

We see that the involution constructed satisfies the requirements of a normal ruling.

\section{Generating families and the homology of the linearized Chekanov-Eliashberg DGA}

Work of Traynor  and collaborators \cite{JTr} \cite{NgTr} \cite{Tr} has used the theory of generating families to distinguish Legendrian knots (actually usually $2$-component links) with identical classical invariants.  A version of their approach is the following:

Let  $F: {\bb R}^N\times {\bb R}\rightarrow {\bb R}$, $f_t = F(\cdot,t)$,  be a linear at infinity generating family for a Legendrian knot $\ell= \ell_F$, and recall from Section 2.5 that $S_F$ denotes the fiber-wise critical set of $F$.
Define $w: {\bb R}^N\times{\bb R}^N\times {\bb R}$ to be the {\sl difference function} 

$$ w(x,y,t) = f_t(x) - f_t(y).$$

Let $\delta >0 $ small enough so that the interval $(0, \delta)$ consists entirely of regular values of $w$.  Define the {\sl generating family homology} of $F$ as the grading shifted homology groups ${\scr G}H_*(F) = H_{*+(N+1)}(w \geq \delta, w = \delta ; {\bb Z}_2)$.  The groups ${\scr G}H_*(F)$ may depend on the choice of generating family for $\ell$.  The invariance statement is analogous to Theorem 2.6.

\begin{proposition}[Jordan-Traynor]  The set of possible (isomorphism types of) graded homology groups $\{ {\scr G}H_*(F) \}$ where $F$ is a linear-quadratic at infinity generating family for $\ell$ is an invariant of Legendrian isotopy.
\end{proposition}

 The more general ``linear-quadratic at infinity'' condition means that outside of a compact subset the generating family $F: {\bb R}^{N_1+N_2} \times {\bb R} \rightarrow {\bb R}$ is a sum, $F(x_1,x_2, t) = L(x_1) + Q(x_2)$ where $L: {\bb R}^{N_1} \rightarrow {\bb R}$ is a non-zero linear function and $Q: {\bb R}^{N_2} \rightarrow {\bb R}$ is a non-degenerate quadratic function.   \cite{JTr} considers the wider class of ``linear-quadratic at infinity'' generating families to allow stabilizations. A {\sl stabilization} increases the dimension of the domain  of a generating family  by summing with a non-degenerate quadratic form in the new variables.

 \paragraph{Remark.} Given a linear-quadratic at infinity generating family $F$ as above, we may obtain a linear at infinity generating family $F'$ with $L_F = L_{F'}$ by precomposing with a diffeomorphism which preserves the fibers ${\bb R}^{N_1+N_2} \times \{t\}$.  Clearly, $ {\scr G}H_*(F)  \cong  {\scr G}H_*(F') $, so either restriction at infinity produces the same set of groups $\{ {\scr G}H_*(F) \}$. \medskip  

  Let ${\ell^s}, 0 \leq s \leq 1 $ be a  Legendrian isotopy and suppose a generating family $F^0: {\bb R}^N\times{\bb R}\rightarrow {\bb R}$ for $\ell^0$ is chosen.  It is shown  in \cite{JTr}, \cite{Tr} that after possibly stabilizing the original generating family  we can find a homotopy $F^s, 0 \leq s \leq 1$ such that $F^s$ is a generating family for $\ell^s$. 
The invariance statement is then proved by showing the generating function homologies $ {\scr G}H_*(F^s), 0 \leq s \leq 1$ are all isomorphic.  The reader is referred to \cite{JTr}, \cite{Tr} for the details of this argument in a slightly different setting.  A key point is that for small enough $\delta>0$ critical values of the corresponding difference functions $w^s$ will not cross into the interval $(0, \delta)$.  

\begin{proposition}

The critical points of the function $w$ with positive critical values correspond to the crossings of the $xy$ diagram of $\ell$.  The critical value is equal to the height of the crossing. Moreover, the indices of the critical points are equal to the degrees of the crossings plus $N+1$.

\end{proposition}

\paragraph{{\bf Proof.}}  Since $$\frac{\partial w}{\partial x}=\frac{\partial f_t(x)}{\partial x},\ \frac{\partial w}{\partial y}=-\frac{\partial f_t(y)}{\partial y}, \frac{\partial w}{\partial t}=\frac{\partial}{\partial t}f_t(x)-\frac{\partial}{\partial t}f_t(y),$$the point $(x,y,t)$ is a critical point of $w$ if and only if (1) $x$ and $y$ are critical points of $f_t$, that is, $(x,t),(y,t)\in S_F$, and (2) the $y$-coordinates of the corresponding points of $L_F$ agree. This happens when $x=y$ in which case the critical value is $0$.  When $x \neq y$,  $(x,y,t)$ corresponds to a crossing in the $xy$-diagram, and $w(x,y,t) = f_t(x) - f_t(y)$ is the difference of the two $z$-coordinates.  $(y,x,t)$ is a second critical point of $w$ corresponding to the same crossing, but $w(x,y,t) = - w(y,x,t)$, so exactly one will have positive critical value.  

$$\ind_{(x,y,t)}w=\ind_yf_t+(N-\ind_xf_t)+e$$where $e=1$, if $t$ is a local maximum of the distance between the strands of the points $(f_t(x),t),(f_t(y),t)$, and $e=0$ if $t$ is a local minimum of this distance. Since the degree of the crossing is $\ind_yf_t-\ind_xf_t+e$, the proposition follows.\smallskip

\paragraph{{\bf Remark.}}  Proposition 5.2 suggests an approach to proving Theorem 5.3.  One could hope to directly  compare the Morse complex for the function $w$  and the linearized complex for an augmentation of the Chekanov-Eliashberg DGA (we need to compare only differentials: the chain spaces are the same by Proposition 5.2). 
The difficulty with this approach in general is the lack of a readily apparent augmentation of $\mathbf{A}(\Gamma)$ where $\Gamma$ is the $xy$-diagram of $\ell_F$ itself.  

 \begin{theorem}  
 {\it If $F$ is a linear at infinity generating family for $\ell$ then there exists a graded augmentation $\varepsilon$ for the Chekanov-Eliashberg DGA of $\ell$ such that ${\scr G}H_*(F) \cong H^\varepsilon_{*}(\ell)$.  That is,
$\{ {\scr G}H_*(F) \}\subset\{ H^\varepsilon_*(\ell) \} $.} 
\end{theorem}

\paragraph{{\bf Outline of proof.}} The proof requires associating an augmentation with a generating family.  According to Theorem 2.6, we have the freedom to work with the Chekanov-Eliashberg DGA of a different representative of the Legendrian isotopy class of $\ell_F$. The analogy between the proofs of Theorem 2.4 and Theorem 2.7 suggests a direct route from $F$  to an augmentation for the DGA arising after applying Ng's resolution procedure in combination with the splash construction as in Section 4.1. In Theorem 2.4, Proposition 2.5 is applied to the Morse complexes of the individual functions $f_t$, whereas in Theorem 2.6 the same proposition is applied to complexes arising from $\varepsilon(X^+_m)$.  One could hope to construct an augmentation directly from a generating family so that (after choosing an appropriate family of metrics ) the $\varepsilon(X^+_m)$ are the matrices for the differentials in the morse complex of $f_{t}$ at appropriate values of $t$.  Such an augmentation is constructed in Section 5.2.1.  The value of $\varepsilon$ on  the remaining  generators from splashes, $X^-_m$, 
 reflects bifurcations in the family of Morse complexes $C(f_t)$.
 
The augmentation and, accordingly, also the differential in the corresponding linearized complex are determined by the fiberwise Morse complexes $C(f_t)$ together with their bifurcation data. Therefore, it is convenient to use a method of computing ${\scr G}H_*(F)= H_*( {\bb R}^N \times{\bb R}^N\times {\bb R}, \{ w \leq \delta\} ; {\bf Z}_2 )$  based on this same information.  Such a method is given in section 5.1.  We replace  $( {\bb R}^N \times{\bb R}^N\times {\bb R}, \{ w \leq \delta\}) $ with the fiberwise descending manifold of the fiberwise critical set $S_w= S_F*S_F$ and collapse the sublevel set $\{ w \leq \delta\}$ to a point.  The resulting space is compact, and we use the fiberwise ascending and descending manifolds of points in $S_F$ to provide it with a CW-complex structure.  This is done in Sections 5.1.3-5.1.4 after recalling the standard bifurcations of Morse complexes from \cite{L} in 5.1.1-5.1.2.  The resulting cellular chain complex has roughly the same number of generators as the linearized complex, and its differential is described in terms of the individual Morse complexes $C(f_t)$ and there bifurcations.  

The remaining sections are dedicated to showing that the two complexes have isomorphic homology groups, and the motivations are purely algebraic.  In 5.1.5 and 5.2.3 we take quotients of both complexes by certain acyclic subcomplexes.  This is done so that the two quotient complexes will have exactly the same size.  After setting up a one-to-one correspondence between generators, we see that the differential on the quotient complexes agree on many generators and disagree in a predictable way on the others.  The proof is completed in Section 5.3 by showing that after taking a further homology preserving quotient the complexes become isomorphic.

\paragraph{Remark.} A close variant of Theorem 5.3. is proved  in \cite{NgTr} and \cite{JTr} for certain classes of $2$-component links where the invariants are explicitly calculated.  

\subsection{Cell decompositions from 1-parameter families} Due to the linear at infinity assumption, for large enough $T$ the spaces $\{w \leq \delta\} \cap \{t \geq T\}$ ( resp. $\{w \leq \delta \} \cap \{t \leq -T\}$ ) is a product $ V \times [T, +\infty)$ (resp. $V \times (-\infty, -T]$) where $V$ is a half-space of ${\bb R}^N\times{\bb R}^N$.  Therefore, for the purpose of computing the homology groups in the statement of Theorem 5.3 we can and will henceforth assume our generating family $\{f_t\}$ is defined on a compact interval $t\in [-T, T]$ with $S_F \subset {\bb R}^N\times(-T, T)$.  

We will need to pair our generating family $\{f_t\}$ with a one parameter family of metrics $\{g_t\}, -T\leq t\leq T$ which we assume to be Euclidean outside of some compact set.   Then we can consider the fiber-wise negative gradient flow $\Phi_s: {\bb R}^N\times [-T,T] \rightarrow {\bb R}^N\times [-T,T]$  generated by the vector field $X_{(x,t)} = (- \nabla_{g_t}f_t)_x$. This flow  preserves the fibers ${\bb R}^N\times{t}$ and by the linearity assumption is globally defined.  For $(x,t)$, if $\lim_{s\rightarrow \pm \infty } \Phi_s(x,t)$ exists it belongs to $(S_F)_t$.  Alternatively, $\Phi_s(x,t)$ eventually follows a line with $\lim_{s\rightarrow \pm \infty} f_t(\Phi_s(x,t)) = \mp \infty$.  

\begin{definition}Given a subset $B \subset S_F$, let ${\scr D}(B) = \{ (x,t) | \lim_{s\rightarrow - \infty } \Phi_s(x,t) \in B\}$ and  ${\scr A}(B) = \{ (x,t) | \lim_{s\rightarrow + \infty } \Phi_s(x,t) \in B\}$ denote the {\it fiber-wise  descending} and {\it ascending} manifolds of $B$.  
\end{definition}

\subsubsection{The Morse complex} Our main tool for computations is an extension of a beautiful perspective on Morse theory originally explored by Ren\'e Thom.  On a compact manifold $M$ suppose we are given a single Morse function $f$ with critical points $\crit(f) = \{p_1, \ldots, p_n \}  $ labeled so that $f(p_1) \geq \ldots \geq f(p_n)$.  Fix a metric $g$ on $M$ satisfying the Morse-Smale condition:

For all $i, j,$ ${\scr D}(p_i) $ and ${\scr A}(p_j)$ intersect transversally.

In this case the descending manifolds will form the cells of a CW-complex structure, $M = \coprod {\scr D}(p_i)$.  The dimension of ${\scr D}(p_i)$ is given by the Morse index $\ind(p_i)$.  
  
If $\ind(p_i)=\ind(p_j)+1$, then the incidence number $\eta(i,j)=[{\scr D}(p_i)\colon{\scr D}(p_j)]$ is the number of descending gradient trajectories from $p_i$ to $p_j$ (which comprise ${\scr D}(p_i)\cap{\scr A}(p_j)$). These trajectories should be counted with signs, but since we consider only modulo 2 homology, we disregard the signs and view $\eta(i,j)$ as integers modulo 2. 

To simplify some of our formulas we will make the convention  that when $\ind(p_i) \neq \ind(p_j) +1$, $\eta(i,j) =0$.

We arrive at the {\it Morse complex} $(C_*(f), \partial)$ with coefficients in ${ \bf Z}_2$.  
$$C_l(f) = {\bf Z}_2 \{p_i | \ind(p_i) = l\}$$
$$ \partial : C_l(f) \rightarrow C_{l-1}(f),  \partial p_i = \sum_{ind(p_j)=l-1} \eta(i,j)p_j $$
Note that the differential is strictly upper triangular in the basis $\{p_i\}$. The homology of the Morse complex is isomorphic to the singular homology of $M$ modulo 2.

If we fix two regular values $a$ and $b>a$ of $f$ and restrict the complex to the critical points $p_i$ with $a<f(p_i)<b$, then we get a complex $(C_\ast(f;a,b),\partial)$ whose homology is $H_\ast(\{f\le b\},\{f\le a\};{\bb Z}_2)$. The latter will work also in the case when $M$ is non-compact, provided the gradient flow is globally defined.  Indeed, in this case a cellular decomposition of the quotient space
\[
N= (\bigcup_{a<f(p_i)<b}{\scr D}(p_i))/\{f \leq a\}
\]
is provided by the descending manifolds ${\scr D}(p_i)$ together with the collapsed set as an extra $0$-cell.   $(C_\ast(f;a,b),\partial)$ agrees with the corresponding cellular chain complex mod the $0$-cell $\{f \leq b\}$ and thus computes the relative homology
\[
\widetilde{H}_\ast(N; {\bb Z}_2) \cong H_\ast(\bigcup_{a<f(p_i)<b}{\scr D}(p_i), \{f \leq a\}; {\bb Z}_2 ) \cong 
H_\ast(\{f\le b\},\{f\le a\} ; {\bb Z}_2 ).
\]

\subsubsection{Bifurcations of the Morse complex} Consider the following two conditions on a 
$C^\infty $ function $f$:

(i) $f$ is Morse;

(ii) critical values of $f$ are distinct;

These conditions are simultaneously satisfied on a generic (open and dense) subset of the set of all functions.

After possibly making a small perturbation of $F$ preserving both ${\scr G}H_*(F)$ and the Legendrian isotopy type of $\ell_F$ we may assume that our family of functions $\{f_t\}, -T \leq t \leq T $,  is generic in the following sense:

(i) and (ii) hold for $f_t$ at all but  a finite number of values $-T< t_1 \leq t_2 \leq \ldots \leq t_M< T$ where precisely one of (i) and (ii) fails for $f_{t_i}$.  

Further, if $f_{t_m}$ fails to be Morse it is due to the presence of a single degenerate critical point $p$ corresponding to either:

\emph{(B)} a birth of a pair of critical points with adjacent indices $\lambda +1$ and $\lambda$, or

\emph{(D)} a death of a pair of critical points with adjacent indices $\lambda +1$ and $\lambda$.
More formally, in a neighborhood of $p$ in $\R^N\times \R$ one can change  coordinates while preserving $t$ so that 
\[
f_t(x_1, \ldots, x_n) = f(p) + x_1^3 + \epsilon (t-t_m) x_1 + Q(x_2, \ldots, x_n)
\]
 where $Q$ is a non-degenerate quadratic form of index $\lambda$ and $\epsilon= -1$ at a birth and $\epsilon = 1$ at a death.

Also, when condition (ii) fails it will be through:

\emph{(TCV)}  transverse intersection of a single pair of critical values.

\paragraph{{\bf Remark.}}  On the front diagram of $\ell_F$ \emph{(B)} and \emph{(D)} points correspond to left and right cusps respectively.  The $t$-values for which two critical points share a common critical value correspond to crossings.  Note that the genericity assumption \emph{(TCV)} guarantees  $\ell_F$ is an \emph{embedded} Legendrian submanifold.

In either case \emph{(B)} or \emph{(D)} there will be a unique degenerate critical point $p$.  ${\scr D}(p)$ (resp. ${\scr A}(p)$) will be a half disk of dimension $\lambda +1$ (resp. $N-\lambda$ where 
$N=\dim M$) so the Morse-Smale transversality condition still makes sense in this setting.

We can  choose the family of metrics $g_t$ so that the Morse-Smale condition fails at only a finite number of $t$-values at which conditions (i) and (ii) are both satisfied \cite{L}.   The manner in which the Morse-Smale condition fails can be assumed to be a standard `handle slide' as described before Lemma 5.2 below.   We use \emph{(M-S)}  to indicate a value of $t$ where the Morse-Smale condition fails for $(f_t, g_t)$.

\begin{definition}

 {\rm A $t$-value such that one of \emph{(B)}, \emph{(D)}, \emph{(TCV)}, or \emph{(M-S)} occurs will be referred to as a {\it singular $t$-value}. }

\end{definition}

\paragraph{{\bf Remark.}}  Throughout the remainder of the proof we will be considering each of these singularity types on a case by case basis.  The reader is encouraged to concentrate on one type of singularity at a time in order to avoid becoming bogged down.

\medskip

Let us examine case by case how the Morse complex changes when we pass a singular $t$-value and in addition extend the Morse complex to the types of singular pairs $(f_t, g_t)$ described above.  We will haphazardly add $t$'s to our previous notation to indicate which function we are considering.  For instance $p_i$ becomes $p_i(t)$,  $\eta(i,j)$ becomes $\eta^t(i,j)$, and the differential of the Morse complex becomes $\partial_t$.

{\it Case} of no singularity:  If the open interval $(a, b)$ contains no singular $t$-values then the Morse complex is stable through out.  We have $\eta^{t'}(i,j) = \eta^{t''}(i,j)$ for all $a< t' , t'' < b$.

In the rest of the cases, for simplicity of notation, we assume that on the interval $(-2,2)$ there is a lone singularity at $t=0$ of the desired type.  The formulas in Lemma 5.1 and 5.2 come from \cite{L}, although the cell-decompositions at singular $t$-values are not discussed in that reference.

{\it Case (B)}:   Assume the newly born critical points are labeled $p_k(1), p_{k+1}(1)$ with the indices 
$\lambda +1, \lambda$.     $C(f_1)$ has two more generators than $C(f_{-1})$. To express the relationship between the two complexes it is convenient to  forego our usual labeling convention by listing $\crit(f_{-1}) = \{ p_1, \ldots, p_{k-1}, p_{k+2}, \ldots \}$ (compare with the labeling convention of Section 4.1). This notation will be retained in future considerations of \emph{(B)} $t$-values.

Differentials are related by the formulas: \[\begin{array} {rl} \partial_1 p_i &= \partial_{-1} p_i + \eta^1(i, k+1)  \partial_1p_k + \eta^1(\partial_{-1}p_i,  p_{k+1}) p_k\\ &= \partial_{-1}p_i + \sum_j \eta^1(i,k+1)\eta^1(k,j) p_j + \sum_j \eta^{-1}(i,j)\eta^1(j, k+ 1) p_k\\ \partial_1 p_k &= p_{k+1} + \sum_{k+1<j} \eta^1(k,j)p_j\\
\partial_1 p_{k+1} &= \sum_j\eta^1(k,j) \partial_{-1}p_j \end{array}\]

\paragraph{{\bf Remark.}}  For $i\neq k, k+1$, unless $\ind(p_i) = \lambda +2$ or $\lambda +1$ we have simply $\partial_1p_i = \partial_{-1}p_i$.  In the use of $\eta^1(\partial_{-1}p_i,  p_{k+1})$ we extend $\eta^1( \cdot , p_{k+1})$ as a linear form.

Let us record some observations.

\begin{lemma}

{\rm(i)} $C(f_{-1}) \cong C(f_1)/{\bf Z}_2\{p_k ,\partial_1p_k \},\ p_i \mapsto [p_i]$ (isomorphism of complexes).

{\rm(ii)} In $C(f_1)/{\bf Z}_2 \{p_k , \partial_1p_k \}$, $$[p_{k+1}] = \sum_{k+1 <j} \eta^1(k, j)[p_j]$$.

{\rm(iii)}  The map $A: C(f_{-1})\oplus {\bf Z}_2\{ p_k, p_{k+1}\} \rightarrow C(f_1)$, $A(p_i)= p_i + \eta^1(i, k+1) p_k$ for $i \neq k, k+1$, $A(p_k)= p_k$, $A(p_{k+1}) = p_{k+1} + \sum_{k+1 < j} \eta^1(k, j) p_j$ is an isomorphism of complexes where $C(f_{-1})\oplus {\bf Z}_2\{ p_k, p_{k+1}\}$ is a direct sum of complexes with differential defined on the second component as $p_k \mapsto p_{k+1} $.

\end{lemma} 

The descending manifolds with respect to $(f_0, g_0)$ will give a CW-complex structure and we define the `Morse complex' $C(f_0)$ to be the cellular chain complex.  The only abnormality is that the descending manifold of the degenerate critical point is a half disc so contributes $2$ cells.  If we label these cells as $p_k, p_{k+1}$ then $C(f_0)$ will be identical to $C(f_1)$.     

{\it Case} (D):  The situation is symmetric to the case (B).

{\it Case} (TCV):  The only thing that changes is the way we should label critical points.  If the crossing occurs between critical points labeled $p_k, p_{k+1}$  then the matrix $(\eta^{-1}(i,j))$  agrees with $(\eta^1(i,j))$ after conjugating by the permutation matrix of the transposition $(k, k+1)$.  The Morse complex at $t=0$ makes sense as usual.  By convention we let the labeling of critical points at $t=0$ agree with the labeling when $t<0$.

{\it Case} (M-S):  The generic way that this will happen is that at some point a gradient flow line will connect two critical points $p_k, p_l, k<l$ with the same index, $\lambda$. This prevents a CW-complex structure by descending manifolds since the closure of ${\scr D}(p_k)$ will intersect  ${\scr D}(p_l)$.  To rectify this problem we divide the cell ${\scr D}(p_l)$ into three cells ${\scr D}(p_l)^-, {\scr D}(p_l)^+, {\scr D}(p_l)^0$ where \[ \begin{array} {rl} {\scr D}(p_l)^0 &= \overline{{\scr D}(p_k)} \cap {\scr D}(p_l)\\
{\scr D}(p_l)^- &= (\overline{{\scr D}(p_k(t), t<0)} \cap {\scr D}(p_l)  )- {\scr D}(p_l)^0\\
{\scr D}(p_l)^+ &= (\overline{{\scr D}(p_k(t), t>0)} \cap {\scr D}(p_l)  )- {\scr D}(p_l)^0 \end{array} \]

Here ${\scr D}(p_l)^0 $ is a $(\lambda -1)$-cell while both ${\scr D}(p_l)^-$ and ${\scr D}(p_l)^+$ are $\lambda$-cells.   

In the Morse complex for $C(f_0)$, which is defined to be the cellular chain complex obtained from this decomposition, we use $p^0_l,p^-_l, p^+_l $ to denote the three pieces of ${\scr D}(p_l)$.   Relationships between the complexes $C(f_i), i= -1,0,1$ are \[ \begin{array} {c} \partial_{-1} p_l =  \partial_1 p_l = \partial_0p^-_l + \partial_0p^+_l\\ \partial_0 p^+_l = p^0_l + x  , \partial_0 p^-_l =
p^0_l + y \end{array} \] where $x, y \in  \span\{ p_i | i \neq l\}$,

$$ \partial_0p_k = \partial_{-1}p_k + \partial _0 p^-_l = \partial_1p_k + \partial_0 p^+_l,$$
and if $\ind(p_i) = \lambda +1$, then
$$ \partial_0p_i = \partial_{-1}p_i + \eta(i,k)p^-_l = \partial_1p_i + \eta(i,k)p^+_l$$
where all appearances of $p_l$ in the last two formulas should be replaced by $p^-_l + p^+_l$.  

We will use some simple consequences.

\begin{lemma} There are isomorphisms of complexes:

{\rm(i)}  $C(f_{-1}) \cong C(f_0)/{\bf Z}_2 \{ p^-_l, 
\partial_0 p^-_l \} $ given by $p_i \mapsto [p_i], i \neq l$ and $p_l \mapsto [p_l^+] $;

{\rm(ii)}  $C(f_{1}) \cong C(f_0)/{\bf Z}_2 \{ p^+_l ,
 \partial_0 p^+_l \} $ given by $p_i \mapsto [p_i], i \neq l$ and $p_l \mapsto [p_l^-] $;

{\rm(iii)}  $A: C(f_{-1}) \stackrel{\cong}{\rightarrow} C(f_1) $ given by $A(p_i) = p_i, i\neq k$ and $A(p_k) = p_k + p_l $.

\end{lemma}

\subsubsection{A cell decomposition for ${\scr D}(S_F)$}

 Let us refine our partitioning of the interval $[-T, T]$ to  $-T= t_0 < t_1< \ldots < t_M = T$ to include all singular $t$-values of the family $(f_t, g_t)$ and one non-singular $t$-value between any two of the singular values.

In  the remainder of Section 5.1 we describe a CW-decomposition of ${\scr D}(S_w)/\{ w \leq \delta\}$ whose corresponding (reduced) cellular chain complex will be proved quasi-isomorphic to a linearized complex of the Chekanov-Eliashberg DGA of a knot Legendrian isotopic to $\ell_F$.  The  cells will be contained either within a single $t$-value, $t_m$, or in the intermediate intervals $(t_{m-1}, t_m)$, and the isomorphism will map such cells to splash generators of the form $x^+_{m;ij}$ and $x^-_{m;ij}$ respectively.  The reader may significantly simplify matters by concentrating on an interval of the form  $[t_{m-1}, t_{m+1}]$ where $t_m$ is a singular $t$-value and $t_{m-1}$ and $t_{m+1}$ are non-singular then working through the remainder of the argument for each singularity type individually.
 
We begin with a cellular decomposition of ${\scr D}(S_F)/\{ F \leq -C \}$ where  $C$ is large enough so that $F( S_F) \subset (-C, C)$.  
 As described in Section 5.1.2, for each $t_m$ we have a CW-decomposition of $({\scr D}(S_F)\cap\{t= t_m\})/\{ F \leq -C \} $  whose cellular chain complex is the Morse complex $C(f_{t_m})$.    We shorten the notation for cells to $p_i(m) = {\scr D}(p_i(t_m))$.  On an interval $(t_{m-1}, t_m)$ the Morse complex $C(f_t)$ is stable and we can form  cells $P_{i}( m) = {\scr D}(p_i(t) , t_{m-1}<t<t_m )$ of dimension $\ind(p_i) +1$.   Finally, let $\alpha$ denote the $0$-cell corresponding to the collapsed set $\{ F \leq -C \}$.

The following theorem incorporates the bifurcation data to give formulas for the differential.

\begin{theorem}   The decomposition$${\scr D}(S_F)/\{ F \leq -C \}=  (\coprod_m (\coprod_ip_i(m)))\coprod(\coprod_m (\coprod_iP_i(m)))\coprod \alpha$$is compatible with a CW-complex structure. The cellular chain complex modulo $\alpha$ is as a vector space  a direct sum $(\oplus_mA(m)) \oplus (\oplus_m B(m))$ where $A(m)$ is spanned by the cells belonging to $M\times\{t_m\}$, and $B(m)$ is spanned by cells belonging to $M\times (t_{m-1}, t_m)$.   

{\rm(a)}  For each $m$, $A(m) = C(f_{t_m})$ is a sub-complex. 

{\rm(b)}  For each $m$, the differential on $B(m)$ is as a sum of three parts,

\[ \begin{array} {c} \partial_{B(m),B(m)}\colon B(m) \rightarrow B(m)\\
\partial_{B(m),A(m-1)}\colon B(m) \rightarrow A(m-1)\\
\partial_{B(m),A(m)}\colon B(m) \rightarrow A(m) \end{array} \]

The first part acts according to  $$\partial_{B(m),B(m)} P_i(m) = \sum_j\eta^t(i,j)P_j(m),$$
the formula being independent of the choice of $t,\ t_{m-1} < t < t_m$.

The part $\partial_{B(m),A(m-1)}: B(m) \rightarrow A(m-1)$  (resp.
$\partial_{B(m),A(m)}: B(m) \rightarrow A(m)$) is defined depending on the nature of the singularity at $t_{m-1}$  (resp.$t_m$).

In the case, when there is no singularity, $\partial_{B(m),A(m-1)}P_i(m) = p_i(m-1)$  (resp. 
$\partial_{B(m),A(m)}P_i(m) = p_i(m)$).

In the case (B), $\partial_{B(m),A(m-1)}P_i(m) = p_i(m-1) $ (resp. $\partial_{B(m),A(m)}P_i(m)= p_i(m) + 
\eta^{t_m}(i,k+1) p_k(m)$.  (The convention that critical points on the interval $(t_{m-1}, t_m)$ are labeled to omit $p_k$ and $p_{k+1}$ is used here.)

In the case (D), $\partial_{B(m),A(m-1)}P_i(m)= p_i(m) + \eta^{t_{m-1}}(i,k+1) p_k(m)$ under the convention that critical points on the interval $(t_{m-1}, t_m)$ are labeled to omit $p_k$ and $p_{k+1}$ (resp. $\partial_{B(m),A(m)}P_i(m)= p_i(m)$).

In the case (TCV), $\partial_{B(m),A(m-1)}P_i(m)= p_{\sigma(i)}(m-1)$ where $\sigma = (k, k+1)$ (resp. $\partial_{B(m),A(m)}P_i(m)= p_i(m)$).  (The asymmetry is due to the labeling convention.)

In the case (M-S):

if $i \neq k,l$, $\partial_{B(m),A(m-1)}P_i(m)= p_i(m-1)$, (resp. $\partial_{B(m),A(m)}P_i(m)= p_i(m)$);

if $i = k$, $\partial_{B(m),A(m-1)}P_k(m)= p_k(m-1) + p^+_l(m-1)$,  (resp. $\partial_{B(m),A(m)}P_k(m)= p_k(m) + p^-_l(m)$);

if $i = l$, $\partial_{B(m),A(m-1)}P_l(m) = p^-_l(m-1) + p^+_l(m-1)$  (resp. $\partial_{B(m),A(m)}P_l(m) = p^-_l(m) + p^+_l(m)$).

\end{theorem}


\subsubsection{Cellular chain complex on the fiber product}

 Let $w_t$ denote the family corresponding to difference function $w:{\bb R}^N\times{\bb R}^N\times{\bb R} \rightarrow {\bb R}$, so that $w_t(x,y) = f_t(x) - f_t(y)$.  Rather than apply  Theorem 5.4 directly to the pair $(w_t, g_t\oplus g_t)$, we observe that $ {\scr D}(S_w)= {\scr D}(S_F)*{\scr D}(S_{-F}) $ and use fiber products of the cells arising from the families $f_t$ and $-f_t$ as described  in Section 5.1.3.  
   
Note that the descending manifolds of $-f_t$ are simply the ascending manifolds of $f_t$.   
As in the previous section we get cell decompositions for both ${\scr D}(S_F)$ and ${\scr A}(S_F) (= {\scr D}(S_{-F}) )$.  The notation for the ascending manifold cells will use $q$'s instead of $p$'s, but will otherwise be identical.  We use $A'(m)$ and $B'(m)$ for the vector spaces spanned by $q_i(m)$ and $Q_i(m)$ respectively.
  Note that for a non-singular $t$-value (and hence for any $t$-value except for those where {\it(M-S)} fails) the Morse complexes $C(f_t)$ and $C(-f_t)$ are dual:  $\partial_tq_i = \sum_j \eta^t(j,i) q_j$.   At birth-death $t$-values the role of $p_k$ (resp. $p_{k+1}$) is played by $q_{k+1}$  (resp. $q_k$).  Similarly, at $t$-values where the Morse-Smale condition fails the role of $p_k$ (resp. $p_l$) is played by $q_l$ (resp. $q_k$).  

There are product cells $p_i(m)\times q_j(m)$ which will have dimension $N + \ind(p_i) - \ind(p_j)$, and since the cells $P_i( m)$ and $Q_j(m)$ are fibered over the interval as $B^{\ind(p_i)} \times  B^1$ and $B^{N-\ind(p_j)} \times  B^1$ the fiber product $ P_i( m) * Q_j(m)$ is a $1+ N + \ind(p_i) - \ind(p_j)$ cell.  Put together, we get a cell decomposition of ${\scr D}(S_F)*{\scr A}(S_F) = {\scr D}(S_w)$.  A CW-decomposition of ${\scr D}(S_w)/\{w \leq \delta\}$ arises from the above cell decomposition. 
This is because on the cells in the descending (resp. ascending) manifold $F$ decreases (resp. increases) as we move away from the critical set.  
Therefore, on the (fiber) product cells $w$ decreases as we move away from the critical points in the center, and the portion of a cell with $w > \delta $ will itself be a cell when non-empty. The collapsing of $\{w \leq\delta\}$ results in a compact space and adds an additional $0$-cell.

The cellular chain complex $({\bf C}, \partial)$ relative to the $0$-cell $\{w \leq \delta\}$ as a vector space can be realized as a quotient of the direct sum $(\oplus_m A(m) \otimes A'(m)) \oplus (\oplus_m B( m) \otimes B'(m))$ by the subspace generated by 

{\rm (i)} $p_i(m)\otimes q_j(m)$ and $P_i(m)\otimes Q_j(m)$ when $i\geq j$ (including some additional cells with superscripts $+, - ,0$ when the Morse-Smale condition fails ).

{\rm (ii)} $p_k(m)\otimes q_{k+1}(m) $ if $t_m$ is a \emph{(B)}, \emph{(D)}, or \emph{(TCV)} singular $t$-value and $p_k(m)$, $p_{k+1}(m)$ are the offending critical points.

We denote the images of $A(m)\otimes A'(m)$ and $B(m)\otimes B'(m)$ in the quotient as ${\bf A}(m)$ and ${\bf B}(m)$, but we will not use a new notation to distinguish between a generator and its coset.   

Our identification with the cellular chain complex is by  $p_i(m)\times q_j(m)  \leftrightarrow p_i(m)\otimes q_j(m)$ and $ P_i( m) * Q_j(m)\leftrightarrow P_i( m) \otimes Q_j(m)$.  

The ${\bf A}(m)$ are subcomplexes with the usual tensor product differential, 

$$\partial_{A(m)}\otimes {\bf 1}_{A'(m)} + {\bf 1}_{A(m)}\otimes \partial_{A'(m)}$$

The differential on ${\bf B}(m)$ is a sum

\[ \begin{array} {rcl} \partial_{B(m), B(m)}\otimes {\bf 1}_{B'(m)} + & {\bf 1}_{B(m)} \otimes \partial_{B'(m), B'(m)}&:{\bf B}(m) \rightarrow {\bf B}(m) \\ +& \partial_{B(m), A(m-1)} \otimes \partial_{B'(m), A'(m-1)} &:{\bf B}(m)\rightarrow {\bf A}(m-1) \\ + & \partial_{B(m), A(m)} \otimes \partial_{B'(m), A'(m)} & :{\bf B}(m)\rightarrow {\bf A}(m)  \end{array} \]

All of the above maps are well defined on the quotient.

\subsubsection{A quasi-isomorphic quotient}  In this section we take a quotient of the complex constructed in the previous section by an acyclic sub-complex ${\bf E}$.  Motivated by linearized complexes coming from the Chekanov-Eliashberg DGA our goal is to make the number of generators in ${\bf A}(m)$ and ${\bf B}(m)$ roughly the same.  We then provide suggestive notation for a basis of the quotient complex and record the formula for the differential.

The sub-complex ${\bf E}$ is the direct sum of sub-complexes of the ${\bf A}(m)$. The intersection ${\bf E}(m):={\bf E}\cap {\bf A}(m)$ is defined depending on the type of singularity at $t_m$.

{\it Case (B)}. Let 
\[
{\bf F}(m) = \span \{ p_k(m) \otimes q_j(m), p_i(m)\otimes q_{k+1}(m)  \mid k+1 <j \quad \mbox{and} \quad   i< k\}.\]  If follows from Lemma 5.1 (since the coefficient of $p_{k+1}$ in $\partial p_k$ is $1$, the coefficient of $q_k$ in $\partial q_{k+1}$ is $1$, and $p_k(m)\otimes q_{k+1}(m)=0$ in ${\bf A}(m)$) that $\partial$ maps ${\bf F}(m)$ isomorphically onto its image so that ${\bf E}(m) := \span {\bf F}(m)\cup\partial {\bf F}(m)$  is an acyclic subcomplex.  

{\it Case (M-S)}. Let 
\[{\bf F}(m) = \span\{ p^-_l\otimes q_j, p_i\otimes q^-_k \mid l<j \quad \mbox{and} \quad i<k\}.\]  It follows from Lemma 5.2 (since the coefficient of $p^0_l$ in $\partial p^-_l$ is $1$, the coefficient of $q^0_k$ in $\partial q^-_k$ is $1$, and $p_l \otimes q_k = 0$ in ${\bf A}(m)$ regardless of superscripts) that $\partial$ maps ${\bf F}(m)$ isomorphically onto its image so that ${\bf E}(m)=\span {\bf F}(m)\cup\partial {\bf F}(m)$  is an acyclic subcomplex. 

{\it All other cases}: ${\bf E}(m) = \{0\}$.

We now use $\overline{{\bf C}} = (\oplus_m \overline{{\bf A}}(m))\oplus (\oplus_m\overline{{\bf B}}(m) )$ for the quotient complex ${\bf C}/{\bf E}$ with inherited direct sum decomposition.  We use $[x]$ to denote the coset of an element $x \in {\bf C}$. 
Since ${\bf E}$ is acyclic the homology groups of $\overline{{\bf C}}$ compute $H_*({\scr D}(S_w), w\leq \delta; {\bb Z}_2) \cong H_*( {\bb R}^N\times{\bb R}^N\times{\bb R}, w\leq \delta;{\bb Z}_2 )$.  
 
To set up the desired isomorphism of homology groups, at this stage we single out a specific basis for $\overline{{\bf C}}$.  It will be a union of bases $\{{\frak x}^+_{m;i,j} \}$ for $\overline{{\bf A}}(m)$ and $\{{\frak x}^-_{m;i,j}\}$ for $\overline{{\bf B}}(m)$. (The reader may notice a similarity between this notation and the notation in Section 4.1. This is done in purpose: the similarly denoted generators will be put into correspondence with each other during the final stage of the proof of Theorem 5.3.)

Not surprisingly the definition of $\{{\frak x}^+_{m;i,j} \}$ depends on the type of singularity at $t_m$.

{\it Case} \emph{(B)}: ${\frak x}^+_{m;i,j} = [p_i(m) \otimes q_j(m)] \in \overline{{\bf A}}(m),  i<j , \{i,j\}\cap \{k,k+1\} = \phi $. Note that from Lemma 3.1
\[ \begin{array} {c} \displaystyle{[p_{k+1}(m)\otimes q_j(m)]  = \sum_l \eta^{t_{m+1}}(k,l) {\frak x}^+_{m;l,j}}, \\ \displaystyle{ [p_i(m)\otimes q_k(m)] = \sum_l \eta^{t_{m+1}}(l,k+1) {\frak x}^+_{m;i, l} .} \end{array} \]

{\it Case} \emph{(M-S)} :
\[ \begin{array} {c} {\frak x}^+_{m;i,j} = [p_i(m) \otimes q_j(m) ] , i <j , i\neq l , j\neq k, \\
{\frak x}^+_{m;l,j} = [p^+_l(m) \otimes q_j(m)], \ {\frak x}^+_{m;i,k} = [p_i(m) \times q_k^+(m)] . \end{array} \]

{\it All other cases}: ${\frak x}^+_{m;i,j} = p_i(m) \otimes q_j(m) \in \overline{{\bf A}}(m),  i<j $  
where if two critical values $p_k, p_{k+1}$ intersect at $t_m$, ${\frak x}^+_{m;k, k+1}$ is not defined  (This is because $w$ is non-positive on ${\scr D}(p_k)\times {\scr A}(p_{k+1})$).

Finally, for each $m$ and $i<j$ define ${\frak x}^-_{m;i,j} = P_i(m) \otimes Q_j(m)$.

\begin{lemma}
The elements ${\frak x}^-_{m;i,j}, {\frak x}^+_{m;i,j} $ as defined above form a basis for $\overline{{\bf C}}$.
\end{lemma}

To conclude this section we record formulas for the differential of $\overline{{\bf C}}$ with respect to the basis from the lemma.

On the subcomplexes  $\overline{{\bf A }}(m)$ if $t_m$ is a non-singular $t$-value then

$$\partial{\frak x}^+_{m;i,j} = \sum_{i <l<j} \eta^{t_m}(i, l) {\frak x}^+_{m;l,j} + \sum_{i <l<j} \eta^{t_m}(l, j) 
{\frak x}^+_{m;i,l} $$.

If $t_m$ is a singular $t$-value then the differential for $\overline{{\bf A}}(m)$ agrees with the differential for 
$\overline{{\bf A}}(m-1) $  except that the term ${\frak x}^+_{m;k,k+1}$ does not exist if the singular value is a crossing or a right cusp.   ($t_{m-1}$ was assumed to be non-singular if $t_m$ is singular.) This is a consequence of Lemmas 5.1 and 5.2 and in fact ${\bf E}$ was chosen specifically so this would be the case.  

$$\partial{\frak x}^-_{m;i,j} = \widehat{\frak x}^+_{m-1;i,j} +  \sum_{i <l<j} \eta^t(i, l) {\frak x}^-_{m;l,j} + \sum_{i <l<j} \eta^t(l, j) {\frak x}^-_{m;i,l} + {\frak x}^+_{m;i,j}$$
where the term $\widehat{\frak x}^+_{m-1;i,j} \in \overline{{\bf A}}(m-1)$  depends on the type of singularity at $t_{m-1}$.

{\it Case of no singularity}:

$$\widehat{\frak x}^+_{m-1;i,j} = {\frak x}^+_{m-1;i,j}$$

{\it Case} \emph{(B)}: Assuming $\{i,j\} \cap \{k,k+1\} = \phi$, $\widehat{\frak x}^+_{m-1;i,j} = {\frak x}^+_{m-1;i,j}$. Also,
\[ \begin{array} {c} \widehat{\frak x}^+_{m-1;k,j} =\widehat{\frak x}^+_{m-1;i,k+1}=\widehat{\frak x}^+_{m-1;k,k+1}= 0\strut \\
\widehat{\frak x}^+_{m-1;k+1,j}  = \displaystyle{\sum_l \eta^m(k,l){\frak x}^+_{m-1;l,j}}\\
\widehat{\frak x}^+_{m-1;i,k} = \displaystyle{\sum_l \eta^m(l,k+1){\frak x}^+_{m-1;i,l}}  \end{array} \]

{\it Case}  \emph{(D)}:
$$\widehat{\frak x}^+_{m-1;i,j} = {\frak x}^+_{m-1;i,j} + \eta^{m-1}(i,k+1){\frak x}^+_{m-1;k,j} + \eta^{m-1}(k,j) {\frak x}^+_{m-1;i,k+1}$$

{\it Case}  \emph{(TCV)}:
$$\widehat{\frak x}^+_{m-1;i,j} = {\frak x}^+_{m-1;\sigma(i),\sigma(j)} $$
where $\sigma$ is the transposition $(k,k+1)$.

{\it Case}  \emph{(M-S)}: Assuming $i\neq k$, $j\neq l$, $\widehat{\frak x}^+_{m-1;i,j} = {\frak x}^+_{m-1;i,j}$. Also,
\[ \begin{array} {c} \widehat{\frak x}^+_{m-1;k,j} = {\frak x}^+_{m-1;k,j} +{\frak x}^+_{m-1;l,j}\\
\widehat{\frak x}^+_{m-1;i,l} = {\frak x}^+_{m-1;i,l} +{\frak x}^+_{m-1;i,k}\\
\widehat{\frak x}^+_{m-1;k,l} = {\frak x}^+_{m-1;k,l} \end{array} \]

\paragraph{{\bf Remark.}}  The fact that the term of $\partial {\frak x}^-_{m;i,j}$ belonging to $\overline{{\bf A}}(m)$ is always ${\frak x}^+_{m;i,j}$ follows from Lemma 5.1 and 5.2 and the choice of ${\bf E}$.  In the coefficient 
$\eta^t(i,l)$ (resp. $\eta^t(l,j)$) of ${\frak x}^-_{m;l,j}$ (resp. ${\frak x}^-_{m;i,l}$) $t$ needs to satisfies $t_{m-1} < t < t_m$, and is independent of the choice since the Morse complex is stable on the interval $(t_{m-1}, t_m)$.   Again, if there is a crossing or right cusp at $t_m$ then ${\frak x}^+_{m;k,k+1} = 0$.    

\subsection{Construction of the augmentation}

Instead of using the Chekanov-Eliashberg DGA as defined by the $xy$-projection of $\ell_F$ itself, we first apply Ng's resolution procedure to the front projection of $\ell_F$ and add a certain number of splashes (see Section 4.1). 

Specifically, recall that we have made a subdivision $-T \leq t_0 < t_1< \ldots < t_M \leq T$ ($T$ should be large enough so that $f_t$ is linear outside this interval) so that 

(i) every singular $t$-value is a $t_i$ and

(ii)  the sequence of $\{ t_m \}$ alternates between singular and non-singular $t$-values.

For $\Lambda_F'$ we add one splash for each of the $t_m$.  If $t_m$ is a non-singular $t$-value  the splash is contained in a small interval about $t_m$.  
If $t_m$ is a singular $t$-value then place the splash in an interval directly to the left of the singular point.  This is not so important for \emph{(M-S)} singularities since they are not reflected by the front diagram of a knot.  What is important is that at crossings or cusps the splash is placed directly to the left of the singularity.  We get a related partitioning  $-T \leq s_0 < s_1 < \ldots < s_M \leq T$ where $s_m$ is a point chosen from the interval containing the $m$-th splashing.  

Denote the $xy$-projection of the resulting knot by $\Gamma_F$.

\subsubsection{An augmentation for ${\bf A}(\Gamma_F)$}

Refer to section 4.1 for notation of generators and formulas for the differential of the corresponding Chekanov-Eliashberg DGA.  We now construct an augmentation $\varepsilon \colon {\bf A}(\Gamma_F) \rightarrow {\bb Z}_2$. 

On the generators $y_m, z_m$ coming from crossing and right cusp inserts set 
$$\varepsilon (y_m ) = 0, \varepsilon(z_m) =0.$$

We define
$$\varepsilon (x^+_{m;ij}) = \eta^{s_m}(i,j)$$
so that $\varepsilon(X^+_m)$ is the matrix of the differential in $C(f_{s_m})$ with respect to the basis $\{p_i(s_m)\}$.  Since $\partial X^+_m = (X^+_m)^2$, it follows that $\varepsilon(\partial X^+_m)= 0$.  Since two critical points whose critical values meet at a crossing (resp. right cusp) cannot be (resp. must be) joined by a gradient trajectory we see also that $\varepsilon (\partial y_m ) = \varepsilon( x^+_{m-1;k,k+1})= 0$ (resp. $\varepsilon (\partial z_m ) = \varepsilon(x^+_{m-1;k,k+1})+1 = 0$).

In general, $\partial X^-_m = X^+_m(I + X^-_m) + (I +X^-_m) \widetilde{X}^+_{m-1}$ (compare with the similar formulas in Section 4.1) where the definition of $\widetilde{X}^+_{m-1}$ depends on the type of insert appearing between the $(m-1)$-th and $m$-th laminated zones.  In all cases $\varepsilon(\widetilde{X}^+_{m-1})$ is already specified, and turns out to be the matrix of the differential in  a complex closely related to $ C(f_{s_{m-1}})$.  The condition $\varepsilon(X^-_m)=0$ will be satisfied provided $I + \varepsilon(X^-_m)$ is the matrix of an isomorphism between the two relevant complexes.  Such isomorphisms are provided in Section 5.1.2.  Please note that due to our (backwards) conventions for matrices of linear maps, compositions of linear maps correspond to matrix products in the reverse order.  

$\varepsilon(X^-_m)$ is defined depending on the type of singularity at $t_{m-1}$.

{\it Case of no singularity}.  In this case $\varepsilon(\widetilde{X}^+_{m-1}) = \varepsilon(X^+_{m-1}) = \varepsilon(X^+_m)$ and we define $\varepsilon(X^-_m) = 0$.

{\it Case} \emph{(B)}. In this case $\varepsilon(\widetilde{X}^+_{m-1}) $ is the matrix of the split extension of $ C(f_{s_{m-1}})$ by ${\bf Z}_2\{p_k \mapsto p_{k+1}\}$,  provided we use or usual convention for labeling critical points.  We define $\varepsilon(X^-_m)$ in such a way that $I + \varepsilon(X^-_m)$ is the matrix of the isomorphism $A$ from Lemma 5.1:
\[ \begin{array} {c} \varepsilon( x^-_{m;i,k} ) = \eta^{s_m}(i,k+1),  \varepsilon( x^-_{m;k+1,j} ) = \eta^{s_m}(k,j)  \\ \varepsilon(x^-_{m;i,j}) = 0\ \mbox{when}\ i \neq k+1\ \mbox{and}\ j \neq k. \end{array} \]

{\it Case} \emph{(D)}. In this case, $\varepsilon( \widetilde{x}^+_{m-1;i,j}) = \eta^{s_{m-1 }}(i,j) + \eta^{s_{m-1 }}(i,k+1) \eta^{s_{m-1 }}(k,j)$, and we observe that $\varepsilon(\widetilde{X}^+_{m-1}) $ is the matrix of the differential in the quotient $C(f_{s_{m-1}})/{\bf Z}_2\{p_k \mapsto \partial p_k \}$ with respect to the basis $\{[p_1],\ldots ,[p_{k-1}],[p_{k+2}], \ldots\}$.  According to the Lemma 5.1 this matrix will be the same as $\varepsilon(X^+_m)$, so we can define $\varepsilon(X^-_m) =0.$ 

{\it Case} \emph{(TCV)}. In this case, $\widetilde{X}^+_{m-1} =  B \widehat{X}^+_{m-1} B^{-1} $
where $B = P_{(k,k+1)} + y_mE_{k+1,k+1}$, and $B^{-1} = P_{(k,k+1)} + y_mE_{k,k}$.  Here $P_{(k,k+1)}$ is the permutation matrix for the transposition $(k,k+1)$, and $E_{k,k}, E_{k+1,k+1}$ are matrices with a single non-zero entry. $ \widehat{X}^+_{m-1} $ is simply $X^+_{m-1}$ with the entry $x^+_{m-1;k,k+1}$ replaced by $0$.   

Since $\eta^{s_{m-1}}(k,k+1) = \eta^{s_{m}}(k,k+1) = 0 $ necessarily due to the crossing,
$$\varepsilon ( B \widehat{X}^+_{m-1} B^{-1} ) = P_{(k,k+1)} \varepsilon(X^+_{m-1}) P_{(k,k+1)} = \varepsilon(X^+_{m}).$$

The last equality follows since the complexes  $C(f_{s_{m-1}})$ and $C(f_{s_{m}})$ differ only by the ordering of generators.  We define $\varepsilon(X^-_m)= 0$.

{\it Case}  \emph{(M-S)}. In this case, $\widetilde{X}^+_{m-1} =  X^+_{m-1} $.  We define $\varepsilon(X^-_m)$ so that $I + \varepsilon(X^-_m)$  is the matrix of the isomorphism $A: C(f_{s_{m-1}}) \stackrel{\cong}{\rightarrow} C(f_{s_{m}}) $ from Lemma 5.2.,
$$\varepsilon(x^-_{m;k,l}) = 1, \varepsilon(x^-_{m;i,j}) = 0, (i,j)\neq(k,l) $$

In the discussion surrounding the definition of $\varepsilon$ we have proved

\begin{lemma}  $\varepsilon$ is an augmentation of ${\bf A}(\Gamma_F)$.

\end{lemma}

\subsubsection{The differential in the linearized complex}
We record formulas for the differential in the associated linearized complex $({\bf A}^\varepsilon, d^\varepsilon)$.

\[ \begin{array} {c} d^\varepsilon y_m = x^+_{m-1;k, k+1} , d^\varepsilon z_m = x^+_{m-1;k, k+1}, \\
d^\varepsilon x^+_{m;i,j} = \sum_l \eta^{s_m}(i,l) x^+_{m;l,j} + \sum_l \eta^{s_m}(l,j) x^+_{m;i,l}. \end{array} \]

The formulas for $d^\varepsilon x^-_{m;i,j}$ depends on the type of singularity at $t_{m-1}$.

{\it Case of no singularity}:
$$d^\varepsilon x^-_{m;i,j} = x^+_{m-1;i,j} + x^+_{m;i,j} +  \sum_l \eta^{s_m}(i,l) x^-_{m;l,j} + \sum_l \eta^{s_{m-1}}(l,j) x^-_{m;i,l}$$
In this case $\eta^{s_m}(i,j) = \eta^{s_{m-1}}(i,j) = \eta^t(i,j)$ for any $t \in [s_{m-1}, s_m]$.

{\it Case} \emph{(B)}: Assuming $\{i,j\}\cap\{k,k+1\} = \phi$,

\[ \begin{array} {c} d^\varepsilon x^-_{m;i,j}= x^+_{m;i,j} + \eta^{s_m}(k,j) x^+_{m;i,k+1} +x^+_{m-1;i, j} + {^\circ x^-} \\
d^\varepsilon x^-_{m;k,j}= x^+_{m;k,j} + \eta^{s_m}(k,j) x^+_{m;k,k+1} +{^\circ x^-} \\
d^\varepsilon x^-_{m;k+1,j}= x^+_{m;k+1,j} + \sum_l\eta^{s_m}(k,l) x^+_{m-1;l,j} + {^\circ x^-} \\
d^\varepsilon x^-_{m;i,k}= x^+_{m;i,k} + \sum_l\eta^{s_m}(l,k+1) x^+_{m;i,l} + {^\circ x^-} \\
d^\varepsilon x^-_{m;i,k+1}= x^+_{m;i,k+1}  + {^\circ x^-} \\
d^\varepsilon x^-_{m;k,k+1}= x^+_{m;k,k+1}  + {^\circ}{x^-} \end{array} \]
where ${^\circ x^-} \in \span\{x^-_{m;i,j} \}$ is a term which will be irrelevant to our argument.

{\it Case} \emph{(D)}:
\[ \begin{array} {rl} d^\varepsilon x^-_{m;i,j}&= x^+_{m;i,j} + \eta^{s_{m-1}}(i,k+1) x^+_{m-1;k,j} + \eta^{s_{m-1}}(k,j) x^+_{m-1;i,k+1} \\ &+ \sum_{i <l<j} \eta^{s_m}(i, l) x^-_{m;l,j} + \sum_{i <l<j} \eta^{s_m}(l, j) x^-_{m;i,l} +\stackrel{\circ}{z} \end{array} \]
where ${^\circ z} \in \span \{z_m\}$ will be irrelevant to the argument.

\paragraph{{\bf Remark.}}  The term $\sum_{i <l<j} \eta^{s_m}(i, l) x^-_{m;l,j}$ comes from linearizing $X^+_m(I+X^-_m)$ while the similar term  $\sum_{i <l<j} \eta^{s_m}(l, j) x^-_{m;i,l}$ comes from linearizing $(I+X^-_m)\widetilde{X}^+_{m-1}$.  The reason that they both use the Morse complex at $s_m$ is that, as mentioned above, $\varepsilon(X^+_m) = \varepsilon(\widetilde{X}^+_{m-1})$.

{\it Case} \emph{(TCV)}:
\[ \begin{array} {rl} \partial x^-_{m;i,j} &= x^+_{m;i,j} + x^+_{m-1;\sigma(i),\sigma(j)} +  \sum_{i <l<j} \eta^{s_m}(i, l) x^-_{m;l,j} \\ &+ \sum_{i <l<j} \eta^{s_m}(l, j) x^-_{m;i,l} +\stackrel{\circ}{y} \end{array} \]
where ${^\circ y}  \in \span\{y_m\}$ will be irrelevant to the argument.

{\it Case}  \emph{(M-S)}. Assuming $i\neq k$, $j\neq l$,
\[ \begin{array} {c} d^\varepsilon x^-_{m;i,j} = x^+_{m;i,j} + x^+_{m-1;i,j} +{^\circ x^-} \\
d^\varepsilon x^-_{m;k,j} = x^+_{m;k,j} + x^+_{m-1;k,j} + x^+_{m-1;l,j} + {^\circ x^-} \\
d^\varepsilon x^-_{m;i,l} = x^+_{m;i,l} + x^+_{m-1;i,l} + x^+_{m;i,k} + {^\circ x^-} \\
d^\varepsilon x^-_{m;k,l} = x^+_{m;k,l} + x^+_{m-1;k,l} + {^\circ x^-} \end{array} \]
where ${^\circ x^-} \in \span\{x^-_{m;i,j} \}$ is a term which will be irrelevant to our argument.

\subsubsection{Another quasi-isomorphic quotient}The extra generators $y_m, z_m$ from crossings and right cusps generate an acyclic subcomplex which has the basis $\{ y_m, x^+_{m-1;k,k+1} \}\cup \{ z_m, x^+_{m-1;k,k+1}\}$ where $m$ ranges over all values such that the $m$-th insert is a crossing or right cusp.  We will work with the quotient and retain our previous notation so that an element now denotes its coset in the quotient.  The remaining basis elements are then exactly the same as the basis elements for the quotient of the cellular chain complex constructed in Section 5.1.5.  

\subsection {The isomorphism between homology groups}

For the remainder of the proof we use the correspondence ${\frak x}^\pm_{m;i,j}\leftrightarrow x^\pm_{m;i,j}$ to view the two complexes as being defined on the same vector space $\overline{{\bf C}}$ (with a grading shift by $N+1$).  We denote the differential inherited from the linearized complex by $d^\varepsilon$ and the differential inherited from the cellular chain complex as $\partial$.  We use the direct sum decomposition $\overline{{\bf C}} = (\oplus_m \overline{{\bf A}}(m)) \oplus (\oplus_m\overline{{\bf B}}( m))$ from section 5.1.5. The proof of Theorem 5.3 will be completed by the following lemma.

\begin{lemma}  The homology groups of $(\overline{{\bf C}}, d^\varepsilon)$ and $(\overline{{\bf C}}, \partial)$ are isomorphic.

\end{lemma}

\paragraph{{\bf Proof of Lemma}} Observe, that  $d^\varepsilon$ and $\partial$ are identical on $\overline{{\bf A}}(m)$ for all $m$, and also on $\overline{{\bf B}}(m)$ in all cases except when $t_{m-1}$ is a singular $t$-value of type \emph{(B)} or \emph{(M-S)}.  

Let 
$$D_1 = \span\overline{{\bf B}}(m)\cup d^\varepsilon (\overline{{\bf B}}(m)),  D_2 = \span\overline{{\bf B}}(m) \cup \partial (\overline{{\bf B}}(m)) $$
where $m$ ranges over all values with $t_{m-1}$ a type \emph{(B)} or \emph{(M-S)} singular $t$-value.

The result will follow from two claims.

Claim 1:  $D_1$ and $D_2$ are acyclic subcomplexes.

Claim 2:  $D_1 = D_2$

\noindent since then the identity map on the quotient will be an isomorphism of complexes. 

Claim 1 is easy to verify since the composition of either differential with the projection to $\overline{{\bf A}}(m)$ is upper triangular with respect to a proper choice of ordering on the bases.  

For Claim 2, let $p_1 = \pi_A\circ d^\varepsilon$ and $p_2= \pi_A\circ \partial$ where 
\[ \pi_A :  \overline{{\bf A}}(m-1) \oplus \overline{{\bf B}}(m) \oplus \overline{{\bf A}}(m) \rightarrow \overline{{\bf A}}(m-1)\oplus \overline{{\bf A}}(m) \]
denotes the projection. Clearly
\[\begin{array}{cc} D_1= \span\{x^-_{m;i,j} , p_1x^-_{m;i,j}\}, & D_2= \span\{x^-_{m;i,j} , p_2x^-_{m;i,j}\} \end{array}\]

To see that $D_1 \subset D_2$ observe that

{\it Case} \emph{(B)}. Assuming $\{i,j\}\cap\{k,k+1\}= \phi$,
\[ \begin{array} {c} p_1x^-_{m;i,j} = p_2 x^-_{m;i,j} +\eta^{s_m}(k,j)p_2 x^-_{m;i,k+1} \\
p_1x^-_{m;k,j} = p_2 x^-_{m;k,j} + \eta^{s_m}(k,j)p_2 x^-_{m;k,k+1}  \\
p_1x^-_{m;k+1,j}= p_2 x^-_{m;k+1,j}  \\
p_1x^-_{m;i,k}= p_2x^-_{m;i,k} + \sum_l\eta^{s_m}(l,k+1) p_2x^-_{m;i,l} \\
p_1x^-_{m;i,k+1}= p_2 x^-_{m;i,k+1} \\
p_1x^-_{m;k,k+1}= p_2 x^-_{m;k,k+1} \end{array} \]

{\it Case} \emph{(M-S)}.  Assuming $j \neq l$
\[\begin{array} {c} p_1x^-_{m;i,j} = p_2x^-_{m;i,j}  \\
p_1x^-_{m;i,l} = p_2x^-_{m;i,l}+ p_2x^-_{m;i,k} \end{array}\]

\paragraph{{ \bf Remark.}}  It should be apparent at this point why we don't explicitly need to know the terms ${^\circ x^-}, {^\circ y}, {^\circ z}$ in the formulas for $d^\varepsilon$.

\section{The Sabloff duality is the Alexander duality}

In \cite{SD} Sabloff established a duality theorem for the linearized contact homology groups $H^\varepsilon_*(\ell)$.

\begin{theorem} [Sabloff]  If $\ell$ is a Legendrian knot and $\varepsilon : {\mathbf{ A}} \rightarrow {\bb Z}_2$ any graded augmentation, then we have

\[\begin{array} {rll} \dim_{{\bb Z}_
2}H^\varepsilon_k(\ell)&= \dim_{{\bb Z}_2}
H^\varepsilon_{-k}(\ell) & k\neq \pm 1 \\ \dim_{{\bb Z}_
2}H^\varepsilon_1(\ell)&= \dim_{{\bb Z}_2}
H^\varepsilon_{-1}(\ell) +1. &  \end{array}\]

\end{theorem}

Together with our Theorem 5.3, this statement gives the following relation for generating family homology groups.

\begin{theorem}

Suppose $F:\R^N\times\R \rightarrow \R$ is a linear at infinity generating family for a Legendrian knot $\ell$.

{\rm (i)} If $k\ne\pm1$, then ${\scr G}H_k(F)\cong {\scr G}H_{-k}(F).$

{\rm (ii)} ${\scr G}H_1(F) \cong {\bb Z}_2 \oplus {\scr G}H_{-1}(F)$.

\end{theorem}

In this section we provide a proof of Part (i) from the generating family perspective. Possibly, Part (ii) can be proven by a careful examination of the homomorphisms in the long exact sequence used below.

(i) is restated in terms of the level sets of the difference function $w$ as follows:

If $k\ne\pm1$, then $H_{N+k+1}(w_{\ge\delta},w_\delta;{\bb Z}_2)\cong H_{N-k+1}(w_{\ge\delta},w_\delta;{\bb Z}_2).$

\medskip

\begin{lemma}  In the region $\delta \geq w \geq -\delta$,  $w$ is a Morse-Bott function.  There is a single non-degenerate critical submanifold $\Delta$ of index $N$.  Furthermore, $\Delta$ is diffeomorphic to $S^1$ and contained in the level $w=0$.

\end{lemma}

\paragraph{Proof.}  It is easy to see that for small enough $\delta$ the critical set will be  $\Delta = \{ (x,x,t) \in {\bb R}^{2N+1} | (x,t) \in S_F \}$.  The Hessian matrix at $(x,x,t) \in \Delta$ has the block form
\[H_{(x,x,t)} = 
\left [  \begin{array} {ccc} 
	A & 0 & b \\ 0 & -A & -b \\ b^T & -b^T & 0
\end{array} \right ]
 \]  where $A = A^T = \left[ \displaystyle \frac{\partial^2 F}{\partial x_i \partial x_j}(x, t)  \right]$ and $b = \left[\displaystyle \frac{\partial^2 F}{\partial x_i \partial t}(x, t) \right]$.  Under the transversality assumption on $F$,  $\begin{array} {rl} [A & b]  \end{array}$ is of full rank. To verify that $\Delta$ is non-degenerate we must check that  $\ker H_{(x,x,t)} = T_{(x,x,t)}\Delta$.
Since \[S_F = \left\{ (x,t) \in {\bb R}^N\times {\bb R} \Big| \frac{\partial F}{\partial x_i}(x,t) = 0, 1 \leq i \leq N \right\},\ T_{(x,t)}S_F = \Ker [A\ b] ,\] \noindent and hence \[T_{(x,x,t)}\Delta =  \{ (\xi, \xi, \tau) \in T_{(x,x,t)}({\bb R}^N\times {\bb R}^N \times  {\bb R}) | (\xi, \tau) \in \Ker[A\ b] \}.\]  

\noindent Now, 

\[ (\xi, \eta, \tau) \in \ker H_{(x,x,t)} \Leftrightarrow \left [ \begin{array}{c} A\xi + \tau b \\ -A\eta -\tau b \\  b^T (\xi -\eta) \end{array} \right ]  = 0 \eqno{6.1}\]

We see immediately that $T_{(x,x,t)}\Delta \subset \ker H_{(x,x,t)}$. For the reverse inclusion, suppose that $(\xi,\eta,\tau )\in \ker H_{(x,x,t)}$.  From (6.1) we see that $(\xi, \tau), (\eta, \tau) \in \ker \begin{array} {rl} [A & b]  \end{array} $.  In addition, $\xi-\eta \in \ker \begin{array} {rl} [A & b]^T  \end{array} = \{0\} \Rightarrow \xi= \eta$, so $T_{(x,x,t)}\Delta \supset \ker H_{(x,x,t)}$ holds.

For the index computation let $T_{(x,x,t)}{\bb R}^{2N+1} = h^+\oplus h^- \oplus \ker H_{(x,x,t)} $ where the direct sum is orthogonal with respect to $H_{(x,x,t)}$ and the hessian is positive (resp. negative) definite when restricted to $h^+$ (resp. $h^-$).   Such a decomposition exists for any symmetric bilinear form defined on a vector space over ${\bb R}$, and $\ind(w, \Delta) = \dim h^-$ is well defined.  Now, note that the isomorphism $S: T_{(x,x,t)}{\bb R}^{2N+1} \rightarrow T_{(x,x,t)}{\bb R}^{2N+1}, S(\xi, \eta, \tau) = (\eta, \xi, \tau) $ satisfies $H_{(x,x,t)}( S u , Sv) = - H_{(x,x,t)}(u, v)$ for any $u, v \in T_{(x,x,t)}{\bb R}^{2N+1}$.  Therefore, $T_{(x,x,t)}{\bb R}^{2N+1} = S(h^-)\oplus S(h^-) \oplus \ker H_{(x,x,t)} $ is another $H_{(x,x,t)}$-orthogonal direct sum, and now $S(h^-)$ is positive definite and $S(h^+)$ is negative definite.  It follows that $\dim h^+ = \dim h^- \Rightarrow \ind(w, \Delta) = N$.

\begin{corollary}  $H_j(w_{\leq\delta}, w_{\leq -\delta}; {\bb Z}_2) \cong H_{j-N}(S^1, {\bb Z}_2)$, hence
\[\dim H_j(w_{\leq\delta}, w_{\leq -\delta}; {\bb Z}_2) = \left \{ \begin{array} {ll} 1, & \mbox{if } j=N,N+1 \\ 0 & \mbox{otherwise}. \end{array}\right. \]

\end{corollary}

\paragraph{Proof.} It follows from Lemma 6.1 and a fundamental result from the Morse-Bott theory that $w_{\leq \delta}$ is homotopy equivalent to $w_{\leq -\delta}$ with the total space of a disk bundle $E_- \rightarrow \Delta$ of dimension $\ind(w, \Delta)$ attached along $\partial E_-$.  We have

$$H_j(w_{\leq \delta}, w_{\leq -\delta}; {\bb Z}_2) \cong H_j(E_-, \partial E_-; {\bb Z}_2) \cong H_{j-N}(S^1; {\bb Z}_2)$$
where the last $\cong$ is the Thom isomorphism.

\begin{lemma}

If $r\ne N-1,N,N+1$, then the inclusion homomorphism $$H_r(w_{\le-\delta};{\bb Z}_2)\to H_r(w_{\le\delta};{\bb Z}_2)$$ is an isomorphism.

\end{lemma}

\paragraph{Proof.} Corollary 6.2 and the homological sequence of the pair $(w_{\le-\delta}, w_{\le\delta})$.

\medskip

The following technical proposition is used to apply Alexander duality in our non-compact setting.

\begin{proposition} \label{prop:tech} There exists a fiber preserving change of coordinates $\varphi: \R^{2N}\times \R \rightarrow \R^{2N}\times \R$, $\varphi(x,y,t) = \left(\Phi(x,y,t), t \right)$, so that $\alpha = w \circ \varphi$ has $\alpha(x,y,t) = x_1 -y_1$ outside of a compact subset. 
\end{proposition}

\noindent {\bf Proof.}
We can assume $F(x,t) = x_1$ outside of a subset $K \subset R^N\times \R$ of the form\footnote{$D(0, r) \subset \R^{2N}$ denotes a closed disk centered at $0$ with radius $r$.} $K = D(0, R_1) \times [-T, T]$ and that $F(x,t) = x_1$ on an open subset\footnote{ This assumption eases smoothness considerations for functions defined in a piecewise manner below.}  containing $\partial K$.  $F(x,t) -x_1$ is compactly supported and hence uniformly bounded, and since $[x_1- y_1] -w(x,y,t)= F(x,t) -x_1 - (F(y,t)-y_1)$ we can find $C>0$ such that
\[
[x_1- y_1] -w(x,y,t) < C \quad \mbox{for all} \,\, (x,y,t) \in \R^{2N +1}.
\]

We will make use of a smooth cut off function $\beta : [0, +\infty) \rightarrow [0,1]$ satisfying
\begin{itemize}
\item $\beta(r) = \left\{\begin{array}{lc} 0,& \mbox{if} \,\, r \leq 2 R_1 \\ 1, & \mbox{if} \,\, r \geq R_2 \end{array} \right. $,

 where  $R_2$ is chosen large enough so that we can also arrange
\item $| \beta'(r)| < \frac{1}{2C}$ for all $r \in [0, +\infty)$.
\end{itemize}

Using $\beta$ to blend $w$ into $x_1-y_1$, we define
\[
\alpha= w +  \beta(r) \left( [ x_1-y_1] - w \right)
\]
where $r(x,y,t) = |(x,y)|$  is the radial coordinate on the fiber.  Then, for $s\in (-\epsilon, 1+\epsilon)$ ($\epsilon >0$ is small) we set 
\[
w_s = (1 -s) w + s \alpha = w + s \beta(r) \left( [x_1 - y_1] - w\right),
\]
and hope to find fiber preserving diffeomorphisms 
\begin{equation}
\label{goal}
\Psi_s : \R^{2N +1} \stackrel{\cong}{\rightarrow} \R^{2N+1} \quad \mbox{with} \quad w_s \circ \Psi_s = w, \quad s \in [0,1].
\end{equation}
Then, we may take $\varphi = (\Psi_1)^{-1}$ to complete the proof.

We will realize the $\Psi_s$ as the flow of a time dependent vector vield $V_s = \frac{d}{ds} \Psi_s$.  Differentiating (\ref{goal}) with respect to $s$ gives
\begin{equation}
\label{vfgoal}
d(w_s)(V_s) +  \frac{d}{d s}(w_s) = 0.
\end{equation}
(\ref{goal}) is satisfied by $\Psi_0 = \mathit{Id}$ when $s=0$, so if we can find $V_s$ satisfying (\ref{vfgoal}) with flow, $\Psi_s$,  defined for $s \in [0,1]$ then (\ref{goal}) may be deduced via an integration.  

$V_s$ will be given in the form $V_s = A_s \frac{\d}{\d x_1} + B_s \frac{\d}{\d y_1}$, hence the $\Psi_s$ will indeed be fiber preserving. (\ref{vfgoal}) then simplifies to
\begin{equation}
\label{newvfgoal}
- \frac{d}{d s}(w_s) = \frac{\d w_s}{\d x_1} A_s + \frac{\d w_s}{\d y_1} B_s
\end{equation}

Now, note that
\[
 \frac{d}{d s}(w_s) = \beta(r) \left( [x_1-y_1] -w \right)
\]
is uniformly bounded and supported\footnote{If $|x|,|y| \leq R_1$ then $\beta(r) = 0$. If $|x|, |y|\geq R_1$ then $[x_1 - y_1] - w =0$.  } in $N_x \cup N_y$  where
\[
N_x = \{ (x,y,t) | \, \, |x| > R_1 , \,  |y|  < R_1 \}, \quad \mbox{and} \quad N_y = \{ (x,y,t) | \, \, |x| < R_1 , \,  |y|  > R_1 \}.
\]
Now, on $N_x$, $w(x,y,t) = x_1 - F(y,t)$, and we compute
\[
\frac{\d w_s}{\d x_1} = 1 + s \beta'(r) \frac{\d r}{\d x_1} \left( [x_1 - y_1] - w\right)
\]
which is uniformly bounded below as
\[
|\frac{\d w_s}{\d x_1}|\geq 1 - (1+\epsilon)| \beta'(r) |C \geq 1/2- \epsilon/2.
\]
Therefore, 
\[
 A_s = \left\{\begin{array}{lc} -\frac{d}{d s}(w_s)/\frac{\d w_s}{\d x_1} & \mbox{on} \,\, N_x 
 \\                             0 & \mbox{elsewhere,} \end{array} \right.
 \]
defines a uniformly bounded smooth function on all of $R^{2N+1}\times (-\epsilon, 1+\epsilon)$.
Similarly, we define
\[
 B_s = \left\{\begin{array}{lc} -\frac{d}{d s}(w_s)/\frac{\d w_s}{\d y_1} & \mbox{on} \,\, N_y 
 \\                             0 & \mbox{elsewhere.} \end{array} \right.
\] 
Clearly, (\ref{newvfgoal}) holds and the flow $\Psi_s$ is globally defined for $s \in [0,1]$ since $V_s$ is uniformly bounded on $\R^{2N+1}\times (-\epsilon, 1+\epsilon) $.

\medskip

\begin{lemma}

$\widetilde H_{N+k}(w_{\ge-\delta};{\bb Z}_2) \cong \widetilde H_{N-k}(w_{\le-\delta};{\bb Z}_2).$


\end{lemma}

\paragraph{Proof:} 

We let $M_+ = w_{\ge-\delta}$ and $M_- = w_{\le-\delta}$.  We consider the $1$-point compactifications $\widehat{M}_+$, $\widehat{M}_+ \subset S^{2N+1}$.  According to Proposition \ref{prop:tech} outside of a compact subset $M_+$ is a closed half space and hence outside of the same subset $\widehat{M}_+$ is a closed ball.  A closed half space deformation retracts to a closed ball (and this can be done leaving an arbitrary compact subset fixed), so $M_+$ and $\widehat{M}_+$  are homotopy equivalent (as are $M_-$ and $\widehat{M}_-)$.  The Alexander duality then gives
\[
\widetilde H_{N+k}(M_+;{\bb Z}_2) \cong \widetilde H_{N+k}(\widehat{M}_+;{\bb Z}_2) \cong \widetilde H_{N-k}(\widehat{M}_-;{\bb Z}_2)
\cong \widetilde H_{N-k}(M_-;{\bb Z}_2).
\]

\medskip



\paragraph{Proof of Theorem 6.3 (i)} If $k=0$, we have nothing to prove. Let $k\ne 0,\pm1$.  

\[ \begin{array} {ll} \hskip12pt H_{N+k+1}(w_{\ge\delta},w_\delta;{\bb Z}_2)\\ \cong H_{N+k+1}({\bb R}^{2N+1},w_{\le\delta};{\bb Z}_2)&\mbox{(excision)}\\ \cong\widetilde H_{N+k}(w_{\le\delta};{\bb Z}_2)&\mbox{(homological\ sequence\ of}\ ({\bb R}^{2N+1},w_{\le\delta})) \\ \cong\widetilde H_{N+k}(w_{\ge-\delta};{\bb Z}_2)&\mbox{(homeomorphism}\ (x,y,t)\mapsto(y,x,t))\\ \cong\widetilde H_{N-k}(w_{\le-\delta};{\bb Z}_2)&\mbox{(Lemma\ 6.3)}\\  \cong\widetilde H_{N-k}(w_{\le\delta};{\bb Z}_2)&\mbox{(Lemma\ 6.2)}\\  \cong H_{N-k+1}({\bb R}^{2N+1},w_{\le\delta};{\bb Z}_2)&\mbox{(homological\ sequence\ of}\ ({\bb R}^{2N+1},w_{\le\delta})) \\ \cong H_{N-k+1}(w_{\ge\delta},w_\delta;{\bb Z}_2)&\mbox{(excision)}\end{array} \]

\end{document}